\input amstex
\magnification=1200
\documentstyle{amsppt}
\vsize 18.8cm

\def\eps{\epsilon}
\def\om{\Omega}
\def\co{\colon}
\def\R{\Bbb R}
\def\d{\roman{d}}
\def\s{s}
\def \R{\Bbb R}
\def\cl{\overline}

\def\N{\Bbb N}
\def\loc{\roman{loc}}
\newcount\fornumber\newcount\artnumber\newcount\tnumber
\newcount\secnumber
\def\Ref#1{%
  \expandafter\ifx\csname mcw#1\endcsname \relax
    \warning{\string\Ref\string{\string#1\string}?}%
    \hbox{$???$}%
  \else \csname mcw#1\endcsname \fi}
\def\Refpage#1{%
  \expandafter\ifx\csname dw#1\endcsname \relax
    \warning{\string\Refpage\string{\string#1\string}?}%
    \hbox{$???$}%
  \else \csname dw#1\endcsname \fi}

\def\warning#1{\immediate\write16{
            -- warning -- #1}}
\def\CrossWord#1#2#3{%
  \def\x{}%
  \def\y{#2}%
  \ifx \x\y \def\z{#3}\else
            \def\z{#2}\fi
  \expandafter\edef\csname mcw#1\endcsname{\z}
\expandafter\edef\csname dw#1\endcsname{#3}}
\def\Tag#1#2{\begingroup
  \edef\mmhead{\string\CrossWord{#1}{#2}}%
  \def\writeref{\write\refout}%
  \expandafter \expandafter \expandafter
  \writeref\expandafter{\mmhead{\the\pageno}}%
\endgroup}

\openin15=\jobname.ref
\ifeof15 \immediate\write16{No file \jobname.ref}%
\else                   \input \jobname.ref \fi \closein15
\newwrite\refout
\openout\refout=\jobname.ref

\def\newsection{\global\advance\secnumber by 1\fornumber=0\tnumber=0}

\def\dff#1..{%
     \global\advance\fornumber by 1
     \Tag{#1}{(\the\secnumber.\the\fornumber)}
 \unskip                   \the\secnumber.\the\fornumber
   \ignorespaces
}%
\def\dffs#1..{%
     \global\advance\fornumber by 1
     \Tag{#1}{\the\secnumber.\the\fornumber}
 \unskip                   \the\secnumber.\the\fornumber
   \ignorespaces
}%

\def\rff#1..{\Ref{#1}}
\def\dft#1..{%
     \global\advance\tnumber by 1
     \Tag{#1}{\the\secnumber.\the\tnumber}
 \unskip                   \the\secnumber.\the\tnumber
   \ignorespaces
}%
\def\rft#1..{\Ref{#1}}

\def\dfa#1..{%
     \global\advance\artnumber by 1
     \Tag{#1}{\the\artnumber}
 \unskip                   \the\artnumber
   \ignorespaces
}%
\def\rfa#1..{\Ref{#1}}

\def\dep#1..{\Tag{#1}{}}
\def\rep#1..{p. \Refpage{#1}}

\topmatter
\title
 Inertial Manifolds on Squeezed Domains
\endtitle
\author Martino Prizzi\\
and\\ Krzysztof P. Rybakowski\endauthor \leftheadtext{Martino
Prizzi and Krzysztof P. Rybakowski} \rightheadtext{Inertial
Manifolds on Squeezed Domains}
\address Universit\`a degli Studi di Trieste, Dipartimento di Scienze Matematiche,
Via Valerio 12/b, 34100 Trieste, Italy
\endaddress \email prizzi\@mathsun1.univ.trieste.it\endemail

\address Universit\"at Rostock, Fachbereich Mathematik, Universit\"atsplatz
1, 18055 Rostock, Germany\endaddress
\email krzysztof.rybakowski\@mathematik.uni-rostock.de\endemail
\date \enddate
\thanks \endthanks

\abstract Let $\om$ be an arbitrary smooth bounded domain in $\R^2$ and $\eps>0$ be
arbitrary.  Squeeze $\om$ by
the factor $\eps$ in the
$y$-direction to obtain the squeezed domain $\om_\eps=\{\,(x,\eps y)\mid
(x,y)\in\om\,\}$.    In this paper we study the family of reaction-diffusion equations
$$ \alignedat 2
u_t&=\Delta u+
f(u),&\quad &t>0,\,
 (x,y)\in\Omega_\eps
\\ \partial _{\nu_\eps} u&=0,&
& t>0,\, (x,y)\in\partial\Omega_\eps,\endalignedat\tag $E_\eps$ $$
where $f$ is a dissipative nonlinearity of polynomial growth.
In a previous paper we showed
 that, as $\eps\to 0$,  the equations $(E_\eps)$ have a limiting
equation which is an
abstract semilinear parabolic equation defined on a closed
linear subspace of
$H^1(\om)$. We also proved  that the family ${\Cal A}_\eps$ of the corresponding attractors
is upper semicontinuous at $\eps=0$. In this paper we prove that, if     $\om$ satisfies some natural
assumptions,  then
 there is a family
$\Cal M_\eps$ of inertial $C^1$-manifolds for $(E_\eps)$ of
some fixed finite dimension $\nu$. Moreover, as $\eps\to 0$, the flow    on $\Cal
M_\eps$ converges in the $C^1$-sense to the limit flow on $\Cal
M_0$.
\endabstract
\endtopmatter

\document\newsection\head \the\secnumber. Introduction\endhead

Let  $\om$ be an arbitrary smooth bounded domain in $\R^2$ and
$\eps>0$ be arbitrary. Write $(x,y)$ for a generic point of
$\R^2$. Given $\eps>0$ squeeze $\om$ by the factor $\eps$ in the
$y$-direction to obtain the squeezed domain $\om_\eps$. More
precisely, let $T_\eps\co \R^2\to \R^2$, $(x,y)\mapsto (x,\eps y)$
and $\om_\eps:=T_\eps (\om)$.

Consider the following reaction-diffusion equation on $\om_\eps$:
$$ \alignedat 2 u_t&=\Delta u+ f(u),&\quad &t>0,\,
 (x,y)\in\Omega_\eps
\\ \partial _{\nu_\eps} u&=0,& & t>0,\,
(x,y)\in\partial\Omega_\eps.\endalignedat\tag $1_{\eps}$ $$
Here, $\nu_\eps$ is the exterior normal vector field on $\partial
\om_\eps$ and $f\co \R\to \R$ is a   $C^1$-nonlinearity of
polynomial growth such that $\limsup_{|s|\to \infty}f(s)/s\le
-\zeta$ for some $\zeta>0$.
 These hypotheses imply  that $(1_\eps)$ generates a  semiflow
$\tilde\pi_\eps=\tilde \pi_ {\eps,f}$ on $H^1(\om_\eps)$ which has
a global attractor $\tilde\Cal A_\eps=\tilde\Cal A_{\eps,f}$.

As $\eps\to 0$ the thin domain $\om_\eps$ degenerates to a
one-dimensional interval.

One may ask what happens in the limit to the family
$(\tilde\pi_\eps)_{\eps>0}$ of semiflows and to the family
$(\tilde\Cal A_\eps)_{\eps>0}$ of attractors. Is there a limit
semiflow and a corresponding limit attractor?

 This problem was first considered by Hale and
Raugel in \cite{\rfa HaRau1..} for the case when  the domain $\om$
is the {\it ordinate set\/} of a smooth positive function $g$
defined on an interval $[a,b]$, i.e. $$\om=\{\, (x,y)\mid
\text{$a<x<b$ and $0<y<g(x)$\,\}.}$$ The authors prove that, in
this case, there exists a limit semiflow $\tilde\pi_0$,  which is
defined by the one-dimensional boundary value problem $$
\alignedat 2 u_t&=(1/g)(gu_x)_x+ f(u),&\quad &t>0,\,
 x\in\left]a,b\right[
\\   u_x&=0,& & t>0,\, x=a,\,b.\endalignedat\tag $1_0$ $$ Moreover,
$\tilde\pi_0$ has a global attractor $\tilde \Cal A_0$ and, in
some sense, the family $(\tilde \Cal A_\eps)_{\eps\ge 0}$ is
upper-semicontinuous at $\eps=0$.

Hale and Raugel also prove that one can modify the nonlinearity
$f$ in such a way that each modified semiflow $\tilde \pi'_\eps$
possesses an
 invariant $C^1$-manifold $\tilde\Cal
M_\eps$ of some fixed dimension $\nu$ which includes the attractor
$\tilde \Cal A_\eps$ of the original semiflow $\tilde \pi_\eps$.
The semiflows   $\tilde \pi_\eps$ and $\tilde \pi'_\eps$ coincide
on the attractor $\tilde \Cal A_\eps$.

Moreover, as $\eps\to 0$,    the reduced
flow on $\tilde\Cal M_\eps$ $C^1$-converges to the reduced flow on
$\tilde\Cal M_0$.

 If
the domain $\om$ is not the ordinate set of some function (e.g. if
$\om$ has holes or different horizontal branches) then $(1_0)$ can
no longer be a limiting equation for $(1_\eps)$. Nevertheless, as
it was proved in \cite{\rfa pr..} the family $\tilde \pi_\eps$
still has a limit semiflow. Moreover, there exists a limit global
attractor and the upper-semicontinuity result continues to hold.

In order to describe the main results of \cite{\rfa pr..}
 we first transfer the family $(1_\eps)$ to boundary value problems on the fixed
domain $\om$. More explicitly, we use the linear isomorphism
$\Phi_\eps\co H^1(\om_\eps)\to H^1(\om)$,   $u\mapsto u\circ
T_\eps$, to transform problem $(1_\eps)$ to the equivalent problem
$$ \alignedat 2 u_t&=u_{xx}+\frac 1{\eps^2}u_{yy}+ f(u),&\quad
&t>0,\,
 (x,y)\in\Omega
\\   u_x\nu_1+\frac 1{\eps^2}u_y\nu_2&=0,& & t>0,\,
(x,y)\in\partial\Omega.\endalignedat\tag $2_\eps$ $$ on  $\om$.
Here, $\nu=(\nu_1, \nu_2)$ is the exterior normal vector field on
$\partial \om$.

Note that equation $(2_\eps)$ can be written in the abstract form
$$\dot u+A_\eps u=\hat f(u)$$ where $\hat f\co H^1(\om)\to
L^2(\om)$ is the Nemitski operator
 generated by the function $f$, and $A_\eps$
is the linear operator defined by $$A_\eps u=-u_{xx}-\frac
1{\eps^2}u_{yy}\in L^2(\om)\quad \text{for $u\in H^2(\om)$ with
$u_x\nu_1+\frac 1{\eps^2}u_y\nu_2 =0$ on $\partial \om$. }$$
Equation $(2_\eps)$ defines a semiflow $\pi_\eps$ on $H^1(\om)$
which is equivalent to $\tilde \pi_\eps$ and has the global
attractor $\Cal A_\eps:=\Phi_\eps(\tilde\Cal A_\eps)$, consisting
of the orbits of all full bounded solutions of $(2_\eps)$.

The operator $A_\eps $ is, in the usual way, induced by the
following bilinear form $$a_\eps (u,v):=\int_\om(u_xv_x+\frac
1{\eps^2}u_yv_y)\,\d x\,\d y, \quad u, v\in H^1(\om).$$ Notice
that, for   every fixed $\eps>0$ and $u\in H^1(\om)$, the formula
$$|u|_\eps=\left(a_\eps(u,u)+|u|^2_{L^2(\om)}\right)^{1/2}$$
 defines a norm on $H^1(\om)$ which is equivalent to $|\cdot|_{H^1(\om)}$.
However, $|u|_\eps\to \infty$ as $\eps\to 0^+$ whenever
$u_y\not=0$ in $L^2(\om)$.

In fact, we see that for $u\in H^1(\om)$ $$\lim_{\eps\to
0^+}a_\eps(u,u)=\cases \int_\om u_x^2\,\d x\,\d y,&\text{if
$u_y=0$}\\+\infty,& \text{otherwise.}\endcases$$ Thus the family
$a_\eps(u,u)$, $\eps>0$, of real numbers has a finite limit (as
$\eps\to 0$) if and only if $u\in H^1_s(\om)$, where we define
$$H^1_s(\om):=\{\,u\in H^1(\om)\mid u_y=0\,\}.$$ This is a closed
linear subspace of $H^1(\om)$.

The corresponding limit bilinear form is given by the formula:
$$a_0 (u,v):=\int_\om u_xv_x\d x\,\d y, \quad u, v\in
H^1_s(\om).$$

The form $a_0$ uniquely determines a densely defined selfadjoint
linear operator $$A_0\co D(A_0)\subset H^1_s(\om)\to L^2_s(\om)$$
by the usual formula $$a_0(u,v)=\langle
A_0u,v\rangle_{L^2(\om)},\quad \text{ for $u\in D(A_0)$ and $v\in
H^1_s(\Omega)$.}$$ Here, $L^2_s(\om)$ is the closure of $H^1_s(\om)$ in the
$L^2$-norm, so $L^2_s(\om)$ is a closed linear subspace of
$L^2(\om)$.

It follows that the Nemitski operator $\hat f$ maps the space
$H^1_s(\om)$ into $L^2_s(\om)$. Consequently the abstract
parabolic equation $$\dot u=-A_0u+\hat f(u)\tag $2_0$ $$ defines a
semiflow $\pi_0$ on the space $H^1_s(\om)$. This is the limit
semiflow of the family $\pi_\eps$. In fact, the following results
are proved in \cite{\rfa pr..}:
\proclaim{Theorem A} Let
$(\eps_n)_{n\in\N}$ be an arbitrary sequence of positive numbers
convergent to zero and $(u_n)_{n\in\N}$ be a sequence in
$L^2(\om)$ converging in the norm of $L^2(\om)$ to some $u_0\in
L^2_s(\om)$. Moreover, let $(t_n)_{n\in\N}$ be an arbitrary
sequence of positive numbers converging to some positive number
$t_0$.

Then
$$\bigl|e^{-t_nA_{\eps_n}}u_n-e^{-t_0A_{0}}u_0\bigr|_{\eps_n}\to
0\quad \text{as $n\to \infty$.}$$ If, in addition,  $u_n\in H^1(\om)$
for every $n\in \N$ and if
 $u_0\in H^1_s(\om)$, then
$$|u_n\pi_{\eps _n} t_n-u_0\pi_0 t_0|_{\eps_n}\to 0\quad \text{as
$n\to \infty$.}$$\endproclaim
 The limit semiflow $\pi_0$ possesses a global attractor
$\Cal A_0$. The upper-semicontinuity result alluded to above reads
as follows: \proclaim{Theorem B} The family of attractors
$\left(\Cal A_\eps\right)_{\eps\in[0,1]}$ is upper-semicontinuous
at $\eps=0$ with respect to the family of norms $|\cdot|_\eps$.

This means that $$\lim_{\eps\to 0^+}\sup_{u\in \Cal
A_\eps}\inf_{v\in \Cal A_0}|u-v|_\eps=0.$$
In particular, there exists an $\eps_1>0$
 and an open bounded set $U$ in $H^1(\om)$ including all the
attractors $\Cal A_\eps$, $\eps\in [0,\eps_1]$. \endproclaim

The definition of the linear operator $A_0$, as given above, is
not very explicit. However, as it is shown in \cite{\rfa pr..},
there is a large class of the so-called {\it nicely decomposed\/}
 domains (see Figure 1) on which $A_0$ can be characterized as a
system of one-dimensional second order linear differential operators,
coupled to each other by certain compatibility and Kirchhoff type
balance conditions. In this case, the abstract limit  equation
$(2_0)$ is equivalent to a parabolic equation on a finite graph (see Figure 2).

Let $(\lambda_\nu)_{\nu\in \N}$ be the nondecreasing sequence of the
eigenvalues  of the limit operator $A_0$ (each of the eigenvalues being repeated
according to its multiplicity).
The first main result of the present paper says that given a nicely
decomposed domain $\om$ satisfying some natural additional
condition it follows
 that the sequence $(\lambda_\nu)_{\nu\in \N}$   satisfies the following gap condition:
$$\limsup_{\nu\to \infty}\frac
{\lambda_{\nu+1}-\lambda_{\nu}}{\lambda_{\nu}{}^{1/2}}>0.\tag\dff
gapi..$$
(Cf.
Theorem~\rft gap...)

The hypotheses of Theorem~\rft gap.. are, in particular, satisfied
in the ordinate set case considered in \cite{\rfa HaRau1..}.

Now the second main theorem in this paper is the Inertial Manifold
 Theorem~\rft myt1... This theorem says that
under hypothesis \rff gapi..
 there is an $\eps_0>0$, $\eps_0\le \eps_1$, and there
exists   a family
 $\Cal M_\eps$, $0\le \eps\le \eps_0$ of (inertial)
 $C^1$-manifolds of some finite dimension $\nu$  such
 that, whenever $0\le \eps\le \eps_0$, then   $\Cal A_\eps\subset \Cal M_\eps$
 and the manifold $\Cal M_\eps$ is locally invariant relative to the
 semiflow $\pi_\eps$ on the neighborhood $U$ of the attractor $\Cal
 A_\eps$. Furthermore, as $\eps\to 0$, the reduced flow on the   manifold $\Cal M_\eps$
 converges in the   $C^1$-sense to the reduced flow on $\Cal M_0$.

In particular, our inertial manifold theorem contains, as a
special case, the inertial manifold theorem of Hale and Raugel and
it even improves the latter. In fact, for each $\eps>0$ small
enough, our inertial manifold $\Cal M_\eps$ is (globally)
invariant with respect to some modified semiflow $\pi''_\eps$
coinciding with the original semiflow $\pi_\eps$ on the
neighborhood $U$ of the attractor $\Cal
 A_\eps$. Thus, close to the attractor, $\Cal M_\eps$ is a locally invariant manifold
for the `true' semiflow $\pi_\eps$.  On the other hand, the
modified semiflow considered in \cite{\rfa HaRau1..} is, in
general, different from the original semiflow on every
$H^1$-neighborhood of the attractor.

In order to prove their result, Hale and Raugel first develop
$H^2$-estimates for the attractors, which in turn imply the
corresponding $L^\infty$-estimates. Then they modify the
nonlinearity $f\co \R\to \R$ so as to obtain a bounded
nonlinearity $\tilde f\co \R\to \R$ which induces a globally
 Lipschitzian Nemitski operator from
$L^2(\om)$ to itself.      Then
they prove the existence of inertial manifolds in the space
$L^2(\om)$ using the `cone-squeezing' technique.
This modification of the {\it function\/} $f$ is the main reason
why the modified semiflow $\tilde\pi'_\eps=\tilde \pi_{\eps,\tilde
f}$ coincides with the original semiflow
$\tilde\pi_\eps=\tilde\pi_{\eps,f}$ only on an $H^2$-neighborhood the attractor.

In the present paper we apply the well-known technique for proving existence of
invariant manifolds based on the contraction principle in an appropriate
 space of   functions of
exponential growth. In some cases, this method can be applied in a direct way
 to obtain existence of inertial manifolds (see e.g. \cite{\rfa chow..}, cf. also
 \cite{\rfa ryb..}).
  Unfortunately, such a direct application of this method
to the present case does not work, since the  gap condition \rff gapi.., which is the best
one possible, does not produce
 gaps which are `big enough' to counterbalance the Lipschitz constant produced by the
  given Nemitski operator.
As a consequence, the nonlinear operator $\Gamma_\eps$ used
to define the inertial manifold $\Cal M_\eps$ is {\bf not} a
contraction with respect to the natural norm $|u|_\eps$.

 At this point we use
  a beautiful idea due to
Brunovsk\'y and Tere\v s\v c\'ak  \cite{\rfa brute..}. This idea
consists in replacing the   norm $|u|_\eps$ by the equivalent
norm $$\|u\|_\eps=L|u|_{L^2}+l|u|_{\eps}.$$ This new norm depends
on two positive constants $L$ and $l$ and one tries to determine
these constants in such a way that the operator $\Gamma_\eps$ is a contraction with respect to
the new norm. That such a choice is possible in our case follows by
an application of the Gagliardo-Nirenberg inequality. Moreover,
instead of modifying
 the nonlinearity $f$ we only need to modify  the corresponding
Nemitski operator $\hat f\co H^1(\om)\to L^2(\om)$, so as to obtain
a globally Lipschitzian nonlinear operator $g\co H^1(\om)\to L^2(\om)$   which
coincides with $\hat f$ on $U$. This is the reason why the
modified semiflow coincides on $U$ with the original semiflow.

The reader is referred to the Reference section for some additional papers on
thin domain problems. More extensive bibliography is contained in
the survey paper \cite{\rfa Rau1..} by G. Raugel.

In this paper we use standard notation, writing $\R$ and $\N$
 to denote the set of reals and positive integers,
 respectively.
 We also  identify
 $\R\times \R$ with $\R^{2}$, writing $(x,y)$ for a generic point of
$\R^{2}$. Finally, we denote by ${\Cal L}^1$ and by ${\Cal L}^2$ the Lebesgue measure in $\R$
and $\R^2$ respectively.

\newsection\head \the\secnumber. Preliminaries
\endhead

In this section we      recall a few definitions and notations
from \cite{\rfa pr..}. We also discuss the existence and
differentiability of Nemitski operators generated by real function.

Let $\om$ be a bounded domain in $\R^2$ with Lipschitz
boundary.

 For $\epsilon>0$ define the symmetric bilinear form
$$
a_\epsilon\colon H^1(\Omega)\times H^1(\Omega)\to\R
$$
by
$$
a_\epsilon(u,v):=\int_{\Omega}\left(\nabla_x u\cdot\nabla_x v+
{{1}\over{\epsilon^2}}\nabla_y u\cdot\nabla_y v\right)\,\d x\,\d y.
$$
Let
$b$ be the restriction of the scalar product $\langle\cdot,\cdot\rangle_{L^2(\om)}$
to $H^1(\Omega)\times H^1(\Omega)$.

Since $\om$ has Lipschitz boundary, $H^1(\om)$ is dense in $L^2(\om)$
 and the inclusion operator $H^1(\om)\hookrightarrow L^2(\om)$ is
compact.
Therefore the pair $(a_\eps,b)$ defines a      self-adjoint linear operator
$A_\eps\co D(A_\eps)\subset L^2(\om)\to L^2(\om)$ with domain
dense in $L^2(\om)$.

Moreover,      there are a nondecreasing sequence of eigenvalues
$(\lambda_{\eps,j})_{j\in\N}$       and a
corresponding $L^2$-complete and $L^2$-orthonormal system
$(w_{\eps,j})_{j\in\N}$ of      eigenvectors of the pair $(a_\eps,b)$ (equivalently, of $A_\eps$).

On $H^1(\Omega)$ define the norm
$$
|u|_\epsilon:=\left(a_\epsilon(u,u)+|u|^2_{L^2(\om)}\right)^{1/2}.
$$
This norm is equivalent to $|\cdot|_{H^1(\om)}$ for every fixed
$\epsilon>0$, but $|u|_\epsilon\to\infty$ as $\epsilon\to 0^+$ whenever
$\nabla_y u\not=0$ in $L^2(\Omega)$.

Define the space
$$ H^1_s(\om)=\{\,u\in H^1(\om)\mid \nabla_yu=0\,\}.$$
Note that $ H^1_s(\om)$ is a closed linear subspace of $H^1(\om)$.

Let
$L^2_s(\om)$ to be the closure of the set $H^1_s(\om)$ in $L^2(\om)$.
It follows that $L^2_s(\om)$ is a Hilbert space under the scalar product of $L^2(\om)$.

Note that $|u|_\eps\equiv|u|_{H^1(\om)}$ for $u\in H^1_s(\om)$.

Let $a_0\colon H^1_\s(\Omega)\times H^1_\s(\Omega)\to\R$ be the "limit"
bilinear form defined by
$$
a_0(u,v):=\int_\Omega\nabla u\cdot\nabla v\,\d x\,\d y=
\int_\Omega\nabla_x u\cdot\nabla_x v\,\d x\,\d y.
$$
Let $ b_0$ be the restriction
 of the scalar product
  $\langle\cdot,\cdot\rangle_{L^2(\Omega)}$ to $H^1_\s(\Omega)\times H^1_\s(\Omega)$.

The pair $(a_0,b_0)$ defines a      self-adjoint linear operator
$A_0\co D(A_0)\subset L^2_s(\om)\to L^2_s(\om)$ with domain
dense in $L^2_s(\om)$.

Furthermore, there are a nondecreasing sequence of eigenvalues
$(\lambda_{0,j})_{j\in\N}$   and a
corresponding $L^2_s$-complete and $L^2$-orthonormal system
$(w_{0,j})_{j\in\N}$ of  eigenvectors of the pair $(a_0,b)$ (equivalently, of $A_0$).

It was proved in \cite{\rfa pr.., Th. 3.3} that
$\lambda_{\eps,j}\to \lambda_{0,j}$ as $\eps\to 0$ for $j\in\N$. Moreover,
if $(\epsilon_n)_{n\in\N}$ is an arbitrary sequence of positive numbers
converging to $0$, then the orthonormal system $(w_{0,j})_{j\in\N}$ can be chosen in such a way
that, up to a subsequence, $|w_{\eps_n,j}-w_{0,j}|_{\epsilon_n}\to 0$  as $n\to\infty$ for
$j\in\N$.

In the sequel we will work with a fixed domain and we will
sometimes write $L^p$ for $L^p(\om)$ and $H^1$ for $H^1(\om)$.

Now let $f\co \R\to \R$ be a given nonlinearity. We will
discuss some conditions on $f$ which guarantee that the Nemitski
operator $\hat f$ defined, for $u\co\om\to \R$, by $\hat f(u)= f\circ u$ restricts to a
$C^1$-operator between certain function spaces.

Let us first recall the following
regularity result for Nemitski operators, which is a slight modification of
\cite{\rfa ambro.., Thms 3.4 and 3.7}.

\proclaim{Theorem~\dft t1.1..} Let $\Omega\subset\R^n$ be open and bounded.
Given $\rho,\sigma\in\R$ with $1\leq\rho<\sigma$ set $\beta:=(\sigma/\rho)-1$
(observe $\beta>0$). Suppose $f\in C^1(\R\to\R)$
satisfies
 $$ |f'(s)|\leq C(1+|s|^\beta)\quad\text{for some $C\in
\left]0,\infty\right[$ and all
$s\in\R$.} \tag\dff f1.1.. $$                 Then the Nemitski
operator $$ \hat f\colon L^\sigma(\Omega)\to L^\rho(\Omega),\quad
u\mapsto f\circ u $$ is well-defined and of class                 $C^1$; the
Fr\'echet differential of $\hat f$ is given by $D\hat
f(u)v=(f'\circ u)\cdot v$ for all $u$ and $v\in
L^\sigma(\Omega)$. Finally, there exists a constant $\hat C\in\left[0,\infty\right[$ such
that the following estimates hold: $$ \align &|\hat
f(u)|_{L^\rho}\leq \hat C(1+|u|_{L^\sigma}^{\beta+1})\\ &\|D\hat
f(u)\|_{\Cal L(L^\sigma,L^\rho)}\leq \hat
C(1+|u|_{L^\sigma}^{\beta})
\endalign
$$ for all $u\in L^\sigma(\Omega)$.
\endproclaim

Theorem~\rft t1.1..  implies the following

\proclaim{Corollary~\dft t1.2..} Let  $\Omega$ be an open, bounded
subset of $\R^2$, with Lipschitz boundary. Let $f\in C^1(\R\to\R)$
satisfy the growth estimate $$ |f'(s)|\leq
C(1+|s|^\beta)\quad\text{for $s\in\R$} $$ where $C$ and
$\beta\in\left[0,\infty\right[$ are arbitrary real constants.
Let $F(y):=\int_0^yf(s)\,\d s$ for
$y\in\R$. Then $f\circ u\in L^2(\om)$ whenever $u\in H^1(\om)$. Moreover,
the Nemitski operator $$\hat f\co H^1(\om)\to L^2(\om),\quad
u\mapsto f\circ u, $$ is well-defined and of class $C^1$ on
$H^1(\om)$; the Fr\'echet differential of $\hat f$ is given by
$D\hat f(u)[v]=(f'\circ u)\cdot v$ for all $u$ and $v\in
H^1(\Omega)$. Let $\hat f|_{H^1_\s}$ be the restriction of $\hat
f$ to $H^1_\s(\Omega)$. Then $\hat f|_{H^1_\s}$ is a $C^1$ map
from $H^1_\s(\Omega)$ to $L^2_\s(\Omega)$, with $D(\hat
f|_{H^1_\s})(u)=(D\hat f(u))|_{H^1_\s}$ for $u\in H^1_\s(\Omega)$.
Furthermore, $F\circ u\in L^1(\om)$ whenever $u\in H^1(\om)$ and
the Nemitski operator $$\hat F\co H^1(\om)\to L^1(\om),\quad
u\mapsto F\circ u, $$ is well-defined and of class $C^1$ on
$H^1(\om)$; the Fr\'echet differential of $\hat F$ is given by
$D\hat F(u)[v]=(f\circ u)\cdot v$ for all $u$ and $v\in
H^1(\Omega)$. Finally, there exists a constant $\hat C$ such that
the following estimates hold: $$ \align &|\hat f(u)|_{L^2}\leq
\hat C(1+|u|_{H^1}^{\beta+1})\\ &\|D\hat f(u)\|_{\Cal
L(H^1,L^2)}\leq \hat C(1+|u|_{H^1}^{\beta})\\ &|\hat
F(u)|_{L^1}\leq \hat C(1+|u|_{H^1}^{\beta+2})\\ &\|D\hat
F(u)\|_{\Cal L(H^1,L^1)}\leq \hat C(1+|u|_{H^1}^{\beta+1})\\
\endalign
$$
\endproclaim
\demo{Proof} Without loss of generality we can assume
that $\beta>0$. Set $\rho=2$
and $\sigma=2(\beta+1)$ (resp.          $\rho=1$
and $\sigma=\beta+2$) in Theorem \rft t1.1... By Sobolev's inequality,
$$ |w|_{L^\sigma}\leq C_\sigma|w|_{H^1}\quad\text{for $w\in
H^1(\Omega)$,} $$ where $C_\sigma$ is a real positive constant. Thus an application
of Theorem~\rft t1.1.. shows that the Nemitski operator $\hat f$ (resp. $\hat F$) is
is well-defined, of class $C^1$ and that the above estimates hold.
 We already know that $\hat f(u)\in
L^2_\s(\Omega)$ whenever $u\in H^1_\s(\Omega)$ (Theorem 5.3 in
\cite{\rfa pr..}). In order to complete the proof, we only need to
show that, if $u$ and $v\in H^1_\s(\Omega)$, then $D\hat
f(u)[v]\in L^2_\s(\Omega)$. This follows immediately from the
formula $$ D\hat f(u)[v]=\lim_{t\to 0}{{\hat f(u+tv)-\hat
f(u)}\over{t}} $$ and from the fact that $L^2_\s(\Omega)$ is
closed in $L^2(\Omega)$. \qed\enddemo

For every $\eps>0$ the linear operator $A_\eps$ is sectorial on $X=L^2(\om)$
and the fractional power space $X^{1/2}$ is equal to $H^1(\om)$.
 Similarly, $A_0$ is sectorial on
$X=L_s^2(\om)$ and in this case $X^{1/2}=H_s^1(\om)$.

Let $g\co H^1(\om)\to L^2(\om)$ be Lipschitzian on the bounded subsets of $H^1(\Omega)$ and such
that $g(H^1_s(\om))\subset L^2_s(\om)$.
  It follows from the results of \cite{\rfa He..} that for every $\eps>0$ there exists a
well-defined local semiflow $\pi_{\eps,g}$ on $H^1(\om)$ generated by the semilinear
parabolic equation
$$\dot u+A_\eps u=g(u).\tag \dff kabse..$$
Moreover, there exists a well-defined local
semiflow $\pi_{0,g}$ on $ H^1_s(\om)$ generated by the semilinear parabolic
equation
$$\dot u+A_0 u=g(u).\tag\dff kabs0..$$
In particular, if a function $f\co\R\to \R$ satisfies the
assumptions of Corollary~\rft t1.2.. then, by this corollary, $g:=\hat f\co H^1(\om)\to
L^2(\om)$ is Lipschitzian on the bounded subsets of $H^1(\Omega)$ and $g(H^1_s(\om))\subset
L^2_s(\om)$, so the local semiflows $\pi_{\eps, \hat f}$,
$\eps\ge0$, are well-defined.

\newsection\head \the\secnumber.
Description of the limit problem\endhead
In our previous paper \cite{\rfa pr..} we introduced the class of
the so-called
nicely decomposed domains. We showed that on nicely decomposed domains both
 the operator $A=A_0$ and
its domain of definition can be  explicitly characterized.

In this section we continue the study of nicely decomposed domains
$\om$
and show, in particular, how to characterize the spaces
$H^1_s(\om)$ and $L^2_s(\om)$ (Proposition~\rft image..).
We also obtain a  compact imbedding result (Proposition~\rft
2.0..) which is important in establishing the gap condition~\rff
myfk0.. in section~3.

For the reader's convenience we will first recall the definition
of a nicely decomposed domains.

We say that an open set
$\om\subset
\R^2$ {\it has connected vertical sections\/} if for every $x\in \R$ the
$x$-section
$\om_x$ is connected. Of course, this section is nonempty if and only if $x\in
P(\om)$, where
$P\co\R\times \R\to\R$,
$(x,y)\mapsto x$ is the projection onto the first component.
Note that, given a nonempty bounded domain $\om$ in
$ \R^2 $,
$J_\om :=P(\om)$ is an open interval in $\R$, that is
$J_\om=\left]a_\om,b_\om\right[$, where $-\infty<a_\om<b_\om<\infty$.

Given $a\in \R$ and $\delta\in \left]0,\infty\right[$ we set
$$I_\delta(a):=\left]a-\delta,a+\delta\right[,
\quad\text{$I_\delta^-(a):=\left]a-\delta,a\right[$ and
$I_\delta^+(a):=\left]a,a+\delta\right[$.}$$

\definition{Definition \dft de1..}
Let $\om$, $\om_1$, $\om_2$ be nonempty bounded domains in $ \R^2$. Set
$a_i:=a_{\om_i}$ and $b_i:=b_{\om_i}$, $i=1$, $2$. Given $c\in\R$
we say that $\om_1$ {\it joins $\om_2$ at $c$ in $\om$\/} if
 the following properties hold:\roster
\item $\cl\om_1\cap
\cl\om_2=\{c\}\times [\beta, \gamma]$ where $\beta=\beta_{\om_1,\om_2}$ and
$\gamma=\gamma_{\om_1,\om_2}$ are some real numbers with $\beta<\gamma$;
\item $c=a_{\om_2}=b_{\om_1}$;
\item $\{c\}\times\left]\beta,
\gamma\right[\subset \om$;
\item whenever
$d\in \left ]\beta, \gamma\right[$, then there is a $\delta=\delta(d)>0$
with the property that
$$I_\delta(d)\subset \left]\beta, \gamma\right[$$ and
$$I_\delta^-(c) \times I_\delta(d)\subset \om_1,\quad
I_\delta^+(c) \times I_\delta(d)\subset \om_2.$$\endroster

We say that $\om_1$ and $\om_2$ {\it join at $c$ in $\om$\/} if
$\om_1$  joins $\om_2$ at $c$ in $\om$ or $\om_2$
joins $\om_1$ at $c$ in $\om$.
\enddefinition

\vbox to 4in{\vss\includegraphics{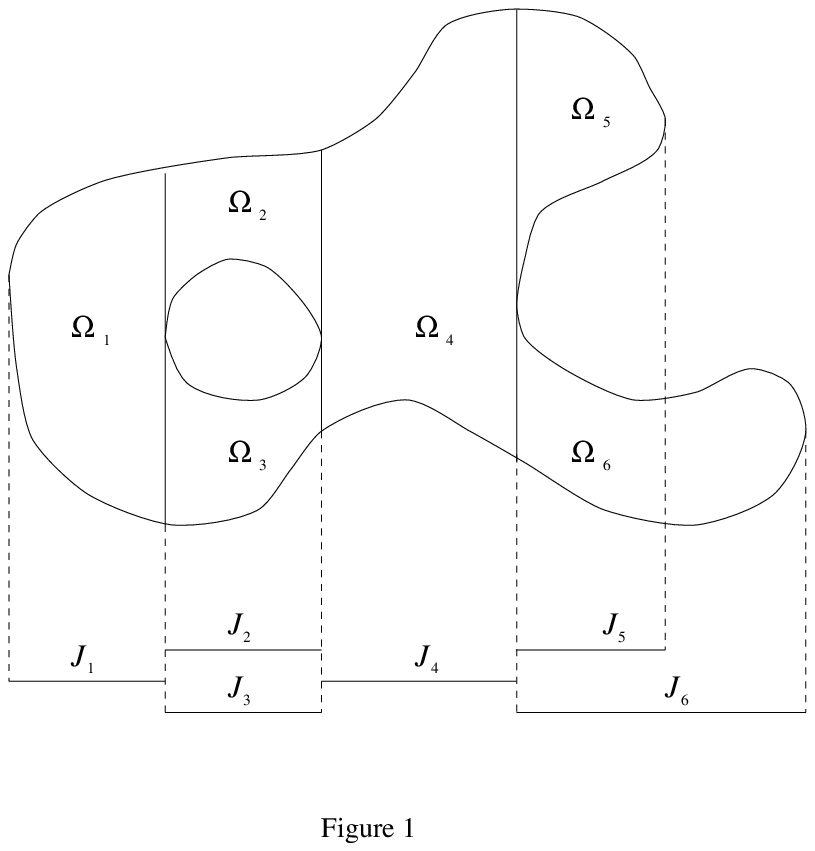}}

\definition{Definition \dft de.., cf Figure 1}
Assume that $\om\subset\R\times \R$ is a nonempty bounded open
domain with Lipschitz boundary.
Let
$P\co \R\times \R\to \R$, $(x,y)\mapsto x$ be the projection onto the first
variable.                  A {\it nice decomposition of\/}
$\om$ is a collection
$\om_1$,
\dots,
$\om_r$ of nonempty pairwise disjoint open connected subsets
of $\om$ with connected vertical
sections                   such that, defining
$J_k:=J_{\om_k}$, $a_k:=a_{\om_k}$, $b_k:=b_{\om_k}$, $k=1$, \dots, $r$,
the following
properties are satisfied:
\roster
\item                  $\om\setminus(\bigcup_{k=1}^r
\om_k)\subset Z$, where $Z:=\bigcup_{l=1}^r(\{a_l,b_l\}\times \R)$;
\item whenever $k=1$, \dots, $r$ then $\partial
\om_k\subset \partial \om\cup (\{a_k,b_k\}\times \R)$ and for
$c\in \{a_k,b_k\}$
$\partial \om_k\cap (\{c\}\times \R)=\{c\}\times I $, where
$I$ is a compact (possibly degenerate) interval in $\R$;
\item                  whenever  $k$, $l=1$, \dots, $r$, $k\not= l$ and
$(c,d)\in\cl\om_k\cap\cl\om_l$ is arbitrary then                   either $\om_k$ and
$\om_l$ join at $c$ in $\om$  or else there is an
$m\in\{1,\dots, r\}$ such that                  $\om_k$ and $\om_m$ join at $c$ in $\om$ and
$\om_l$ and $\om_m$ join at $c$ in $\om$;
\item for every $k=1$, \dots, $r$ the function $p_k\co J_k\to \left]0,
\infty\right[$,
$x\mapsto {\Cal L}^1((\om_k)_x)$, is such that $1/p_k\in L^1(J_k)$.
\endroster\enddefinition

\remark{Remarks}\roster
\item Definition \rft de.. says that,                  up to a set of measure zero,
contained in a set $Z$ of finitely many vertical lines, $\om$ can be
decomposed into the finitely many domains $\om_k$, $k=1$, \dots,
$r$ in such a way that at $Z$ the various sets $\om_k$ and $\om_l$
`join' in a nice way. Points of $\cl\om\cap Z$ are, intuitively speaking,
those at which                  connected components of the vertical sections
$\om_x$ bifurcate.
\item Let $R_1$, \dots, $R_s$ be closed bounded rectangles in $\R^2$
with edges parallel to the coordinate axes, and
$\Omega_R$ be the interior of the union $\bigcup_{k=1}^s R_k$.
 Then, clearly, every connected component of $\Omega_R$
  is a nicely decomposable domain with
 Lipschitz boundary.
\item Let $\phi\colon\R^2\to\R$ be a $C^2$ function, such that
$\phi(x,y)\to\infty$ as $|(x,y)|\to\infty$. Assume that $0$ is
a regular value of $\phi$. Moreover, suppose that, whenever
$\phi(x_0,y_0)=0$ and $\phi_y(x_0,y_0)=0$ then $\phi_{yy}(x_0,y_0)\not=0$.
Let $\Omega_\phi:=\{\,(x,y)\mid \phi(x,y)<0\,\}$ and let $\Omega$
be any connected component of $\Omega_\phi$. It is a nice exercise
to check that $\Omega$ is nicely decomposable.
\item It can be shown that all real analytic domains are nicely decomposable.
\endroster
\endremark

Finally, given a nice decomposition
$\om_1$, \dots, $\om_r$
of $\om$, we set
$$E:=\bigcup_{k=1}\bigl((\{a_k,b_k\}\times \R)\cap\partial \om_k\bigr).$$

Our goal is to give a detailed description of the spaces
$H^1_s(\Omega)$ and $L^2_s(\om)$ when $\Omega$ is a nicely
decomposed domain. We begin our description by first considering
the simpler case of an open set $O$ with connected vertical
sections. Below, such a role will be played by the sets $\Omega_k$
occurring in the nice decomposition of $\Omega$. We do not assume
that $O$ has Lipschitz boundary, since the sets $\Omega_k$
occurring in the nice decomposition of $\Omega$ in general do not
have this property.

Let $P\co \R\times\R\to\R$, $(x,y)\mapsto x$ be the projection onto
the first variable. Let $J:=P(O)$
and assume for simplicity that $J=\left]0,1\right[$.

The following proposition is proved in \cite{\rfa pr..}:

\proclaim{Proposition \dft co..}
Suppose $O$ has connected vertical sections.
Let $J:=P(O)$ and
define the function $p\co J\to \left]0,\infty\right[$ by
$x\mapsto {\Cal L}^1(\om_x)$. If
$u\in L^2(O)$ satisfies $u_y=0$ in the distributional
sense, then there is a set $S$ of measure $0$ in $\R^2$ and a function
$v\in L^1_{\loc}(J)$ such that
$u(x,y)=v(x)$ for every $(x,y)\in O\setminus S$. Moreover,
$p^{1/2}v\in L^2(J)$. If $u\in H^1(O)$ then
$v' \in L^1_{\loc}(J)$ and
 we can choose the set $S$ so
that $u(x,y)=v(x)$ and $u_x(x,y)=v'(x)$
for every $(x,y)\in \om\setminus S$. Moreover,
$p^{1/2}v'\in L^2(J)$
and we can choose the function $v$ to be absolutely continuous on
$J$. The function
$\tilde u\co O\to \R$, $\tilde u(x,y)=v(x)$ is then a continuous
representative of $u$.\qed
\endproclaim

Now, since $O$ is open and bounded, it is easy to check that the function $p$
satisfies the following property:\medskip
(A)
$p\in L^\infty(0,1)$ and  for every
$\epsilon$, $0<\epsilon<1-\epsilon$, there exists $\delta>0$
such that $p(x)>\delta$ a.e.
in $\left]\epsilon,1-\epsilon\right[$.\medskip
Now given an arbitrary function $p$ satisfying property $(A)$ note
 that $p(x)>0$ a.e. in $\left]0,1\right[$.
Therefore the   linear spaces
$$
H(p):=\left\{\,u\in L^1_{\roman{loc}}(0,1)\mid
p^{1/2}u\in L^2(0,1)
\,\right\}
$$
and
$$
V(p):=\left\{\,u\in L^1_{\roman{loc}}(0,1)\mid
u'\in L^1_{\roman{loc}}(0,1), p^{1/2}u\in L^2(0,1),
p^{1/2}u'\in L^2(0,1)\,\right\}
$$
are well-defined.
 Define on $H(p)$ and $V(p)$ the scalar products
$$
\langle u,v\rangle_{H(p)}:=\int_0^1p(x)u(x)v(x)\,\d x
$$
and
$$
\langle u,v\rangle_{V(p)}:=
\int_0^1p(x)u'(x)v'(x)\,\d x+\int_0^1p(x)u(x)v(x)\,\d x.
$$
It is easy to check that these products define Hilbert space
 structures on $H(p)$ and $V(p)$.

Now, given $O$ and $p$ as in Proposition~\rft co.., define the mapping
$$
\iota\co L^2_s(O)\to H(p),\quad u\mapsto v,
$$
where $v$ is the function given by Proposition \rff co... It turns out that
$\iota$ is an isometry of $L^2_s(O)$ onto $H(p)$. Moreover,
$\iota$ restricts to an isometry of $H^1_s(O)$ onto $V(p)$.

We have the following result:

\proclaim{Proposition \dft 2.0..}
Assume that the function  $p$ satisfies (A) and $(1/p)\in L^1(0,1)$.
Let $u\in V(p)$. Then there exists a function
$v\in C^0([0,1])$ such that $u=v$ a.e. in $\left]0,1\right[$.
Moreover, the imbedding $V(p)\hookrightarrow C^0([0,1])$
(and hence the imbedding $V(p)\hookrightarrow H(p)$) is compact.
\endproclaim
\demo{Proof}
Let $\epsilon$, $0<\epsilon<1-\epsilon$, be arbitrary. Then there
exists $\delta>0$ such that $p(x)\geq\delta$ a.e. in
$\left]\epsilon,1-\epsilon\right[$.
It follows that $u|_{\left]\epsilon,1-\epsilon\right[}\in H^1(\epsilon,1-\epsilon)$.
Since $\epsilon$ is arbitrary, we obtain that
there exists a function $v\in C^0(\left]0,1\right[)$ such that
$u=v$ a.e. in $\left]0,1\right[$. We show that $v$
can be extended to a continuous function on $[0,1]$. In fact,
if $x$ and $x'\in\left]0,1\right[$, $x<x'$, then, by H\"older's inequality,
$$
|v(x')-v(x)|^2=|\int_x^{x'}p^{-1/2}p^{1/2}v'\,\d s|^2\leq
\int_x^{x'}(1/p)\,\d s\int_x^{x'}pv'\null^2\,\d s.
$$
It follows that
$$
|v(x')-v(x)|\leq |v|_{V(p)}\left(\int_x^{x'}(1/p)\,\d s\right)^{1/2}.
\tag\dff f2.000..
$$
Since $(1/p)\in L^1(0,1)$, then
$|v(x')-v(x)|\to 0$ as $x,x'\to 0$ (resp. $x,x'\to 1$), so there
exist the limits
$v(0):=\lim_{x\to0}v(x)$ and $v(1):=\lim_{x\to0}v(x)$. By
continuity, estimate \rff f2.000.. holds for $x$ and $x'\in[0,1]$.
Now fix $\cl\epsilon$, $0<\cl\epsilon<1-\cl\epsilon$. Since
$v|_{\left]\cl\epsilon,1-\cl\epsilon\right[}\in H^1(\cl\epsilon,1-\cl\epsilon)$,
it follows that there exists a constant $K_1$ such that
$$
|v|_{L^\infty(\cl\epsilon,1-\cl\epsilon)}
\leq K_1|v|_{H^1(\cl\epsilon,1-\cl\epsilon)}.
$$
Then we can find another constant $K_2$, such that
$$
|v|_{L^\infty(\cl\epsilon,1-\cl\epsilon)}
\leq K_2|v|_{V(p)}.\tag\dff f2.001..
$$
Now \rff f2.000.. and \rff f2.001.. together,
imply that there is a third constant $K_3$ such that
$$
|v|_{L^\infty(0,1)}\leq K_3|v|_{V(p)},\tag\dff f2.002..
$$
so $V(p)\hookrightarrow C^0([0,1])$ with continuous imbedding.
Finally, \rff f2.002.. and \rff f2.000.., and a
straightforward application of Ascoli-Arzela theorem, imply
that the imbedding is compact.
\qed\enddemo

Now we consider the full nicely decomposed domain $\Omega$.
\proclaim{Lemma \dft restr..}
Let $\Omega$ be a nicely decomposed domain. Then, for $k=1$. \dots, $r$,
the following properties hold:
\roster
\item whenever $u\in L^2_s(\Omega)$, then $u|_{\Omega_k}\in L^2_s(\Omega_k)$;
\item whenever $u\in H^1_s(\Omega)$, then $u|_{\Omega_k}\in H^1_s(\Omega_k)$.
\endroster
\endproclaim
\demo{Proof}
Part (2) is obvious and part (1) follows directly from part (2) and from
the definition of $L^2_s(\Omega)$ and $L^2_s(\Omega_k)$.
\qed\enddemo

For $k=1$,\dots, $r$ let us define the linear spaces
$$
H_k:=\left\{\,u\in L^1_{\roman{loc}}(a_k,b_k)\mid
p_k^{1/2}u\in L^2(a_k,b_k)
\,\right\}
$$
and
$$
V_k:=\left\{\,u\in H_k\mid
u'\in L^1_{\roman{loc}}(a_k,b_k),~
p_k^{1/2}u'\in L^2(a_k,b_k)\,\right\}.
$$
We have seen that the spaces $H_k$ and $V_k$ endowed with with the scalar products
$$
\langle u,v\rangle_{H_k}:=\int_{a_k}^{b_k}p_k(x)u(x)v(x)\,\d x
$$
and
$$
\langle u,v\rangle_{V_k}:=
\int_{a_k}^{b_k}p_k(x)u'(x)v'(x)\,\d x+
\int_{a_k}^{b_k}p_k(x)u(x)v(x)\,\d x
$$
respectively, are Hilbert spaces and that the imbedding
$V_k\hookrightarrow H_k$ is dense and compact.
Moreover, consider the following bilinear forms on $V_k$:
$$
a_{\langle k\rangle}(u,v):=\int_{a_k}^{b_k}p_k(x)u'(x)v'(x)\,\d x,
$$
$$
b_{\langle k\rangle}(u,v):=\int_{a_k}^{b_k}p_k(x)u(x)v(x)\,\d x.
$$

Define the product spaces
$$
H_\oplus:=H_1\oplus\cdots\oplus H_r:=\left\{\,
[u]=(u_1,\dots,u_r)\mid u_k\in H_k,k=1,\dots,r
\,\right\}
$$
and
$$
V_\oplus:=V_1\oplus\cdots\oplus V_r:=\left\{\,
[u]=(u_1,\dots,u_r)\mid u_k\in V_k,k=1,\dots,r
\,\right\},
$$
with the scalar products
$$
\langle[u],[v]\rangle_{H_\oplus}:=\sum_{k=1}^r\langle u_k,v_k\rangle_{H_k}
$$
and
$$
\langle[u],[v]\rangle_{V_\oplus}:=\sum_{k=1}^r\langle u_k,v_k\rangle_{V_k}
$$
respectively. It is easy to check that
$H_\oplus$ and $V_\oplus$ are Hilbert spaces and that the imbedding
$V_\oplus\hookrightarrow H_\oplus$ is dense and compact.

Furthermore, consider the following bilinear forms
on $V_\oplus$:
$$
a_\oplus([u],[v]):=\sum_{k=1}^ra_{\langle k\rangle}(u_k,v_k)
$$
and
$$
b_\oplus([u],[v]):=\sum_{k=1}^rb_{\langle k\rangle}( u_k,v_k).
$$
Note that $b_\oplus$ is just the restriction to $V_\oplus\times V_\oplus$
 of the scalar product
$\langle\cdot, \cdot\rangle_{H_\oplus}$.

For $k=1$, \dots, $r$, let us define the mapping $$ \iota_k\co
L^2_s(\Omega_k)\to H_k,\quad u\mapsto v, $$ where $v$ is the
function given by Proposition \rff co... It turns out that
$\iota_k$ is an isometry of $L^2_s(\Omega_k)$ onto $H_k$ and that
$\iota_k$ restricts to an isometry of $H^1_s(\Omega_k)$ onto $V_k$.
Moreover, let us define $$ \iota_\oplus\co L^2_s(\Omega)\to
H_\oplus,\quad \iota_\oplus
u:=(\iota_1(u|_{\Omega_1}),\dots,\iota_r(u|_{\Omega_r})). $$ It follows
 that $\iota_\oplus$ is an isometry of $L^2_s(\Omega)$ into
$H_\oplus$ and that $\iota_\oplus$ restricts to an isometry of
$H^1_s(\Omega)$ into $V_\oplus$.

 Define $$ \gather
V\tilde{_\oplus}:=\\ \biggl\{\, [u]\in V_\oplus\mid
u_k(b_k)=u_l(a_l)~ \text{whenever $b_k=a_l=c$ and $\Omega_k$ and
$\Omega_l$ join at $c$.} \,\biggr\}\endgather $$ Now we are able
to characterize the spaces $H^1_s(\om)$ and $L^2_s(\om)$:

\proclaim{Proposition \dft image..}
The following properties hold:
\roster \item $\iota_\oplus(L^2_s(\Omega))=H_\oplus$; \item
$\iota_\oplus(H^1_s(\Omega))=V\tilde{_\oplus}$. \endroster
\endproclaim
\demo{Proof}
This result is more or less implicitly contained in the proofs of Theorems
6.5 and 6.6 in \cite{\rfa pr..}.
For the reader's convenience, we summarize the main steps of the proof.

First we prove (2). Let $u\in
H^1_s(\Omega)$. We shall prove that $\iota_\oplus u\in
V\tilde{_\oplus}$. Assume that for some $k,l\in\{1,\dots,r\}$,
$b_k=a_l=c$ and $\Omega_k$ and $\Omega_l$ join at $c$. We have to
prove that $$ \iota_k(u|_{\Omega_k})(c)=\iota_l(u|_{\Omega_l})(c).
$$ Set $$ v_k:=\iota_k(u|_{\Omega_k}),\quad
v_l:=\iota_l(u|_{\Omega_l}). $$ Then, by Lemma 6.5 in \cite{\rfa pr..}, we
immediately obtain that $v_k(c)=v_l(c)$, and we are done.

Now assume, conversely, that $[u]\in V\tilde{_\oplus}$. We shall prove that
there exists $v\in H^1_s(\Omega)$ such that $\iota_\oplus v=[u]$.
We define a function $v$ on $\Omega$ in the following way:
$$
v(x,y):=\cases
u_k(x)&\text{for all $(x,y)\in\Omega_k$, $k=1$, \dots, $r$}\\
0&\text{for all $(x,y)\in \Omega\setminus(\bigcup_{k=1}^r\Omega_k)$.}
\endcases
$$
As in  the proof of Theorem 6.6 in \cite{\rfa pr..}, one can show that $v\in
H^1_s(\Omega)$ and that
$$
\alignedat2
v_x(x,y)&=u'_k(x)&\quad&\text{a.e. in $\Omega_k$, $k=1$, \dots, $r$,}\\
v_y(x,y)&=0&\quad&\text{a.e. in $\Omega$.}\\
\endalignedat
$$
Finally observe that, by construction,
$\iota_\oplus v=[u]$ and the proof of (2) is complete.

Now we prove (1). Let $[u]\in H_\oplus$. We shall prove that there exists
$v\in L^2_s(\Omega)$ such that $[u]=\iota_\oplus v$. We define a
function $v$ on $\Omega$ in the following way:
$$
v(x,y):=\cases
u_k(x)&\text{for all $(x,y)\in\Omega_k$, $k=1$, \dots, $r$}\\
0&\text{for all $(x,y)\in \Omega\setminus(\bigcup_{k=1}^r\Omega_k)$}
\endcases
$$
Then $v\in L^2(\Omega)$. We claim that $v\in L^2_s(\Omega)$. This means
that $v$ can be approximated in the $L^2(\Omega)$ norm by functions
of $H^1_s(\Omega)$. Let $\epsilon>0$ be arbitrary. For $k=1$, \dots, $r$,
take $\phi_k\in C^\infty_0(\left]a_k,b_k\right[)$ such that
$$
\int_{a_k}^{b_k}p_k|u_k-\phi_k|^2\,\d x<{\epsilon\over k}.
$$
Define a function $\psi$ on $\Omega$ in the following way:
$$
\psi(x,y):=\cases
\phi_k(x)&\text{for all $(x,y)\in\Omega_k$, $k=1$, \dots, $r$}\\
0&\text{for all $(x,y)\in \Omega\setminus(\bigcup_{k=1}^r\Omega_k)$}
\endcases
$$
By the proof of (2), $\psi\in H^1_s(\Omega)$. Moreover,
$$
\int_\Omega|v-\psi|^2\,\d x\,\d y=\sum_{k=1}^r\int_{\Omega_k}
|v-\psi|^2\,\d x\,\d y=
\sum_{k=1}^r\int_{a_k}^{b_k}p_k|u_k-\phi_k|^2\,\d x<\epsilon,
$$
and our claim is proved. Finally observe that, by definition,
$\iota_\oplus v=[u]$ and the proof is complete.
\qed\enddemo
We end this section with some remarks concerning the self-adjoint
operator  $A$   generated by $a$ in $L^2_s(\Omega)$.
Let $a\tilde{_\oplus}$ be the restriction of $a_\oplus$ to
$V\tilde{_\oplus}\times V\tilde{_\oplus}$
 and let $A\tilde{_\oplus}$ be
 the self-adjoint
operator  generated by
$a\tilde{_\oplus}$ in $H_\oplus$.
  If $u\in D(A)$, then, for all
$v\in H^1_s(\Omega)$, $$ \langle Au,v\rangle_{L^2_s(\Omega)}
=a(u,v)=a\tilde{_\oplus}(\iota_\oplus u,\iota_\oplus v). $$ On the
other hand, $$ \langle Au,v\rangle_{L^2_s(\Omega)}= \langle
\iota_\oplus Au,\iota_\oplus v\rangle_{H_\oplus}. $$ It follows
that $$ a\tilde{_\oplus}(\iota_\oplus u,\iota_\oplus v) =\langle
\iota_\oplus Au,\iota_\oplus v\rangle_{H_\oplus} $$ for all $v\in
H^1_s(\Omega)$, so $\iota_\oplus u\in D(A\tilde{_\oplus})$ and
$A\tilde{_\oplus}\iota_\oplus u=\iota_\oplus Au$. Similarly, one
can prove that, whenever $[u]\in D(A\tilde{_\oplus})$, then
$\iota_\oplus^{-1}[u]\in D(A)$, and $A\iota_\oplus^{-1}[u]=
\iota_\oplus^{-1}A\tilde{_\oplus}[u]$. This means that
$\iota_\oplus$ restricts to an isometry of $D(A)$ onto
$D(A\tilde{_\oplus})$ and that
$A=\iota_\oplus^{-1}A\tilde{_\oplus}\iota_\oplus$.

For $k=1$, \dots, $r$, let us define the spaces $$ Z_k:=\left\{\,
u\in V_k\mid(p_ku')'\in L^1_{\roman{loc}}(a_k,b_k),~
p_k^{-1/2}(p_ku')'\in L^2(a_k,b_k) \,\right\}, $$and let
 $$ Z_\oplus:=Z_1\oplus\cdots\oplus Z_r. $$ It is
not difficult to show that, whenever $u\in Z_k$, then $p_ku'$ can
be extended to a continuous function on $[a_k,b_k]$ (cf. the proof
of Proposition \rft 2.0..). Now we can restate Theorem 6.6 in
\cite{\rfa pr..} in the following way: \proclaim{Theorem \dft
restate..} Let $A\tilde{_\oplus}$ be the self-adjoint operator
generated by the bilinear form $a\tilde{_\oplus}$. Then
$D(A\tilde{_\oplus})=Z\tilde{_\oplus}$, where $Z\tilde{_\oplus}$
is the subspace of $Z_\oplus$ consisting of all
$[u]=(u_1,\dots,u_k)$ satisfying the following properties: \roster
\item $u_k(b_k)=u_l(a_l)$
whenever $b_k=a_l=c$ and $\Omega_k$ and $\Omega_l$ join at $c$;
\item whenever                  $\Gamma$ is a connected component of $E$ (necessarily of
the form $\Gamma=\{c\}\times I$, where $c\in\bigcup_{k=1}^r
\{a_k,b_k\}$ and $I$ is an interval) then $$ \sum_{k\in
\sigma_+}(p_k{u_k}')(c)=\sum_{k\in \sigma_-}(p_k{u_k}')(c). $$
Here, $\sigma_+=\sigma_+(\Gamma)$                  is the set of all $k$ such that
$\cl\om_k\cap (\{b_k\}\times \R)\subset \Gamma$ and so $b_k=c$,
while $\sigma_-=\sigma_-(\Gamma)$ is the set of all $k$ such that
$\cl\om_k\cap (\{a_k\}\times \R)\subset \Gamma$ and so $a_k=c$.
\endroster
Moreover, for $[u]\in Z\tilde{_\oplus}$,
$A\tilde{_\oplus}[u]=(p_1^{-1}(p_1u_r')',\dots,p_r^{-1}(p_ru_r')')$.\qed
\endproclaim
Let $f\colon\R\to\R$ be a $C^1$-function, satisfying the
assumptions of Corollary \rft t1.2... As in Section~2 consider the
abstract semilinear parabolic equation      $$ \dot
u=Au+\hat f(u)\tag\dff abstract.. $$ on $H^1_s(\om)$. Assume that $\Omega$ is
nicely decomposable. Then, due to the isometry $\iota_\oplus$
and in view of Theorem \rft restate..,                  the abstract equation \rff
abstract.. is equivalent to the following system of `concrete'
one-dimensional reaction-diffusion equations: $$
\partial_tu_k=(1/p_k)(p_k{u_k}')'+ f(u_k) \quad\text{on
$\left]a_k,b_k\right[$, for $k=1$, \dots, $r$, } $$ with
compatibility conditions $$ u_k(c)=u_l(c) $$ whenever $b_k=a_l=c$
and $\Omega_k$ and $\Omega_l$ join at $c$, and Kirchhoff type
balance conditions $$ \sum_{k\in
\sigma_+(\Gamma)}(p_k{u_k}')(c)=\sum_{k\in
\sigma_-(\Gamma)}(p_k{u_k}')(c) $$ whenever $\Gamma=\{c\}\times I$
is a connected component of $E$.

\vbox to 2in{\vss\includegraphics{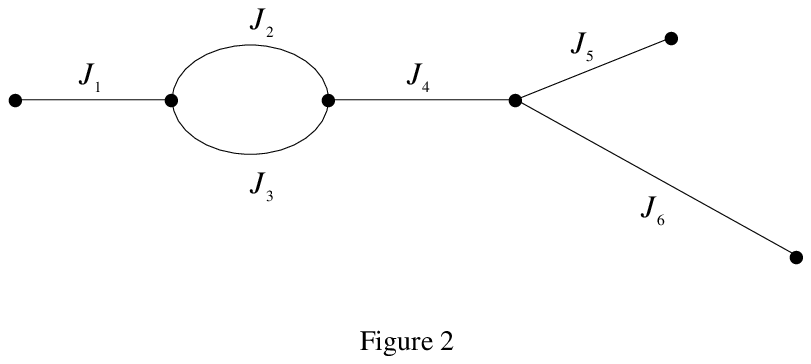}}

As it was  explained in \cite{\rfa pr..}, such a system can be interpreted
as a reaction-diffusion equation on an appropriate finite graph $\Cal G$.
Figure 1 suggests that the edges of this graph are the
intervals $[a_k,b_k]$, for $k=1$,
\dots, $r$  and their endpoints are its vertices.
Moreover, each interval should
be repeated the number of times it occurs in the sequence
$([a_k,b_k])_{k=1}^r$.
On the other hand each endpoint $c$ should be
repeated the number of times it
occurs as the $x$-component of a connected component of $E$ (see Figure 2).

\newsection\head \the\secnumber. Spectral Gap Condition\endhead
In this section
let $\Omega\subset\R^2$ be a nonempty bounded open set with Lipschitz boundary.
We denote by $(\lambda_\nu)_{\nu\in\N}$
 the repeated sequence of the eigenvalues of
the pair $(a,b)$ where
 $$
a(u,v):=\int_\Omega\nabla u\cdot\nabla v\,\d x\,\d y= \int_\Omega
u_x v_x\,\d x\,\d y \quad u,v\in H^1_\s(\Omega) $$ and
 $$b(u,v)=          \int_\Omega uv\,\d x\,\d y
\quad u,v\in H^1_\s(\Omega). $$

We will prove in this section               that, under some additional
assumptions on $\Omega$ (cf. condition (C) below), the following
so-called {\it gap condition\/} holds:
$$\limsup_{\nu\to \infty}\frac
{\lambda_{\nu+1}-\lambda_{\nu}}{\lambda_{\nu}{}^{1/2}}>0,\tag\dff
myfk0...$$
This gap condition will enable us, in the next section,
to establish the                existence of   inertial manifolds
for the semiflows $\pi_\eps$, for $\eps\ge0$ small enough.

We first prove a simple general estimate on the sequence
$(\lambda_\nu)$:
\proclaim{Proposition \dft mypr-1..} There is a $\beta\in
\left]0,\infty\right[$ such that
$$\lambda_\nu\le \beta^2\nu^2 \quad \text{for all
$\nu\in\N$.}$$ \endproclaim
\demo{Proof} There are nonempty and bounded open intervals $I$, $I'$, $J$ and
$J'\subset\R$ such that $$I\times J\subset \om\subset I'\times
J'.$$ Let $\gamma$, $\gamma'$ and $\delta$ be the lengths of $J$, $J'$ and $I$,
respectively. For $\nu\in\N$ let $\Cal V_\nu$ (resp. $\tilde\Cal
V_\nu$) be the set of all $\nu$-dimensional subspaces of
$H^1_s(\om)$ (resp. $H^1_0(I)$).
For every
function $v\in H^1_0(I)$ consider its trivial extension $\tilde v\co I'\to \R$,
$\tilde v(x):=v(x)$ for $x\in I$ and $\tilde v(x):=0$ for $x\in I'\setminus I$. It follows that
$\tilde v\in H^1_0(I')$. Therefore the function
 $u\co I'\times J'\to \R$ defined            by
$u(x,y)=\tilde v(x)$ for $(x,y)\in I'\times J'$ lies in $H^1(I'\times
J')$, so $u|_\om\in H^1(\om)$. By Theorem ~2.5 in \cite{\rfa pr..}
we have $u\in H^1_s(\om)$. The operator $$\Phi\co H^1_0(I)\to
H^1_s(\om),\quad v\mapsto u,$$ is linear and injective. In particular,
$\Phi(\tilde \Cal V_\nu)\subset \Cal V_\nu$, $\nu\in\N$. Thus, for every $\nu\in\N$,
$$\lambda_\nu=\inf_{E\in\Cal V_\nu}\sup_{u\in E\setminus \{0\}}\frac
{a(u,u)}{b(u,u)}\le \inf_{E\in\Phi(\tilde\Cal V_\nu)}\sup_{u\in E\setminus \{0\}}\frac
{a(u,u)}{b(u,u)}$$ $$=\inf_{E\in\tilde\Cal V_\nu}\sup_{v\in E\setminus \{0\}}\frac
{a(\Phi v,\Phi v)}{b(\Phi v,\Phi v)}\le
\inf_{E\in\tilde\Cal V_\nu}\sup_{v\in E\setminus \{0\}}\frac
{\gamma'\int_I v'\cdot v'\,\d x}{\gamma\int_Iv\cdot v\,\d x}=(\gamma'/\gamma)\mu_\nu,$$
where $\mu_\nu$ are the eigenvalues of the Dirichlet problem
$$ \cases v''(x)=-\mu
v(x),&x\in I\\ v(x)=0,&x\in\partial I.\\
\endcases$$
It follows that $\mu_\nu=(\pi/\delta)^2\nu^2$, $\nu\in\N$. The
proposition is proved.\qed
\enddemo
For the rest of this section assume that $\om$
admits a nice decomposition
$\Omega_1$,\dots, $\Omega_r$.           Using  the notation of the last
section we obtain from Proposition~2.2 in \cite{\rfa pr..}
 an eigenvalue-eigenvector sequence $(\lambda_\nu^k,u_\nu^k)_{\nu\in\N}$
 of the
 pair $(a_k,b_k)$ of bilinear forms such that $$
0=\lambda_1^k\le\lambda_2^k\leq\lambda_3^k\leq\dots $$ and
$(u_\nu^k)_{\nu\in\N}$ is a complete orthonormal system in $H_k$,
$k=1$, \dots, $r$.
We also obtain an eigenvalue-eigenvector
 sequence $(\lambda_\nu^\oplus,[u]_\nu^\oplus)_{\nu\in\N}$ of
the pair $(a_\oplus,b_\oplus)$ such that $$
0\leq\lambda_1^\oplus\leq\lambda_2^\oplus\leq\lambda_3^\oplus\leq\dots
$$ and $([u]_\nu^\oplus)_{\nu\in\N}$ is a complete orthonormal system
in $H_\oplus$.

The following result holds:
\proclaim{Proposition \dft mypr-2..}
 $$
\lambda_\nu\geq\lambda^\oplus_\nu\quad\text{for all $\nu\in\N$.}\tag\dff
ungl.. $$\endproclaim
\demo{Proof} The proof is very similar to that of Proposition~\rft mypr-1...
  For $\nu\in\N$ let $\Cal V_\nu$ (resp. $\tilde\Cal
V_\nu$) be the set of all $\nu$-dimensional subspaces of
$V_\oplus$
 (resp. $H^1_s(\om)$).
 Since the map $\Phi:=\iota_\oplus|_{H^1_s(\om)}\co H^1_s(\om)\to V_\oplus$
is an isometry (hence injective) it follows that
$\Phi(\tilde \Cal V_\nu)\subset \Cal V_\nu$, $\nu\in\N$.
Moreover, $a_\oplus(\Phi u,\Phi v)=a(u,v)$ and
$b_\oplus(\Phi u,\Phi v)=b(u,v)$ for $u$, $v\in H^1_s(\om)$.  Thus, for every $\nu\in\N$,
$$\lambda^\oplus_\nu=\inf_{E\in\Cal V_\nu}\sup_{[u]\in E\setminus \{0\}}\frac
{a_\oplus([u],[u])}{b_\oplus([u],[u])}
\le \inf_{E\in\Phi(\tilde\Cal V_\nu)}\sup_{[u]\in E\setminus \{0\}}\frac
{a_\oplus([u],[u])}{b_\oplus([u],[u])}$$ $$=\inf_{\tilde E\in\tilde\Cal V_\nu}\sup_{u\in
\tilde E\setminus \{0\}}\frac
{a_\oplus(\Phi u,\Phi u)}{b_\oplus(\Phi u,\Phi u)}=
\inf_{\tilde E\in\tilde\Cal V_\nu}\sup_{u\in
\tilde E\setminus \{0\}}\frac
{a( u, u)}{b( u,u)}=\lambda_\nu.$$
 This completes the proof.\qed
 \enddemo

Now let $p$ be a function satisfying the  assumptions of Proposition~\rft 2.0...
Consider  the following bilinear form           on $V(p)$: $$
a^p(u,v):=\int_0^1p(x)u'(x)v'(x)\,\d x. $$ Furthermore, define
$b^p(u,v):=\langle u,v\rangle_{H(p)}$ for $u,v\in H(p)$. Then Proposition~\rft
2.0.. and
Proposition~2.2 in \cite{\rfa pr..} implies that there exists a sequence
$(\lambda_\nu^p,u_\nu^p)_{\nu\in\N}$ of eigenvalue-eigenvector pairs of
$(a^p,b^p)$ such that $$
0=\lambda_1^p\le\lambda_2^p\leq\lambda_3^p\leq\dots $$ and
$(u_\nu^p)_{\nu\in\N}$ is a complete orthonormal system in $H(p)$.

 The next proposition describes the
asymptotic behavior of the sequence $(\lambda_\nu^p)_{\nu\in\N}$.

\proclaim{Proposition \dft t2.1..} Let $p\in L^\infty(0,1)$
satisfy one of the following conditions: \roster
\item there exist two constants $\alpha$ and $\beta$, $0<\alpha\leq\beta$,
such that $\alpha \leq p(x)\leq \beta$ almost everywhere in
$\left]0,1\right[$;
\item there exists a function $q\in C^0([0,1])\cap C^2(\left]0,1\right])$,
with $q(x)>0$, $q'(x)\geq0$ and $q''(x)\leq0$ on $\left]0,1\right]$, $q(0)=0$
and $(1/q)\in L^1(0,1)$, and there exist two constants $\alpha$
and $\beta$, $0<\alpha\leq\beta$, such that $\alpha q(x)\leq
p(x)\leq \beta q(x)$ almost everywhere in $\left]0,1\right[$.
\endroster
Let $(\lambda_\nu^p)_{\nu\in\N}$ be the repeated sequence of the
eigenvalues of $(a^p,b^p)$ above. Then there exists a constant
$\gamma>0$, which depends only on $\alpha$ and $\beta$, such that
$$ \lambda_\nu^p\geq\gamma (\nu-1)^2,\quad \nu=1,2,\dots $$
\endproclaim

The proof of Theorem \rft t2.1.. is based on the following

\proclaim{Lemma \dft t2.2..} Let $q\in C^0([0,1])\cap C^2(\left]0,1\right])$,
with $q(x)>0$, $q'(x)\geq0$ and $q''(x)\leq0$ on $\left]0,1\right]$, $q(0)=0$
and $(1/q)\in L^1(0,1)$. Let $(\lambda_\nu^q)_{\nu\in\N}$ be the
repeated sequence of the eigenvalues of $(a_q,b_q)$ above. Then $$
\lambda_\nu^q\geq \pi^2(\nu-1)^2,\quad \nu=1,2,\dots $$
\endproclaim
\demo{Proof} See the Appendix.\qed\enddemo

\demo{Proof of Proposition \rft t2.1..} For every $\nu\in\N$,
$$ \lambda^p_\nu=\inf_{E\in\Cal V_\nu(p)}\sup_{u\in E\setminus\{0\}}
{{a_p(u,u)}\over{b_p(u,u)}}, $$ where $\Cal
V_\nu(p)$ is the set of all $\nu$-dimensional subspaces of $V(p)$. If
condition (1) holds, then $H(p)=L^2(0,1)=H(1)$ and
$V(p)=H^1(0,1)=V(1)$ with equivalent norms. Then $$ \align
\lambda^p_\nu&=\inf_{E\in\Cal V_\nu(p)}\sup_{u\in E\setminus\{0\}}
{{\int_0^1pu'\null^2\,\d x}\over{\int_0^1pu^2\,\d x}}\\
&\geq\inf_{E\in\Cal V_\nu(1)}\sup_{u\in E\setminus\{0\}}
{{\alpha\int_0^1u'\null^2\,\d x}\over{\beta\int_0^1u^2\,\d x}}\\
&=(\alpha/\beta)\lambda^1_\nu,\quad \nu\in\N
\endalign
$$ Now we observe that $$ \lambda^1_\nu=\inf_{E\in\Cal
V_\nu(1)}\sup_{u\in E\setminus\{0\}} {{\int_0^1u'\null^2\,\d
x}\over{\int_0^1u^2\,\d x}},\quad \nu=1,2,\dots $$ are the
eigenvalues of the Neumann problem $$ \cases u''(x)=\lambda
u(x),&x\in\left]0,1\right[\\ u'(x)=0,&x\in\{0,1\}\\
\endcases
$$ This implies that $\lambda^1_\nu=\pi^2(\nu-1)^2$, $\nu=1,2,\dots$,
and the conclusion follows. Assume now that condition (2) holds.
Then $H(p)=H(q)$ and $V(p)=V(q)$ with equivalent norms. Then $$
\align \lambda^p_\nu&=\inf_{E\in\Cal V_\nu(p)}\sup_{u\in
E\setminus\{0\}} {{\int_0^1pu'\null^2\,\d x}\over{\int_0^1pu^2\,\d
x}}\\ &\geq\inf_{E\in\Cal V_\nu(q)}\sup_{u\in E\setminus\{0\}}
{{\alpha\int_0^1 qu'\null^2\,\d x}\over{\beta\int_0^1qu^2\,\d
x}}\\ &=(\alpha/\beta)\lambda^q_\nu,\quad \nu=1,2,\dots
\endalign
$$ By Lemma \rft t2.2.., $\lambda^q_\nu\geq\pi^2(\nu-1)^2$,
$\nu=1,2,\dots$, and the theorem is proved. \qed\enddemo

For $k=1$,\dots, $r$ and $\tau\in\R$, $\tau>0$, let us define $$
m(k,\tau):=\max\left\{\,\nu\in\N\mid\lambda^k_\nu<\tau\,\right\}, $$
and let $$
m(\oplus,\tau):=\max\left\{\,\nu\in\N\mid\lambda^\oplus_\nu<\tau\,\right\}.
$$ We have the following \proclaim{Proposition \dft mult..} Let
$\tau\in\R$, $\tau>0$. Then
$m(\oplus,\tau)=\sum_{k=1}^rm(k,\tau)$.
\endproclaim
\demo{Proof} For $k=1$, \dots, $r$ and $\nu\in\N$, define $$
[u]^k_\nu:=(\underbrace{0,\dots,0}_{k-1},u^k_\nu,0,\dots,0). $$ It is
easy to check that $(\lambda^k_\nu,[u]^k_\nu)$ is an
eigenvalue-eigenvector pair of $(a_\oplus,b_\oplus)$ and that
$([u]^k_\nu)^{k=1,\dots,r}_{\nu\in\N}$ is a complete orthonormal
system in $H_\oplus$. The conclusion follows. \qed\enddemo

Assume now that $\Omega$ satisfies the following additional
hypothesis:
\medskip
(C) For every $k=1$, \dots, $r$ one of the following conditions is
satisfied: \roster
\item there exist two constants $\alpha_k$ and $\beta_k$,
$0<\alpha_k\leq\beta_k$, such that $\alpha_k \leq p_k(x)\leq
\beta_k$ in $\left]a_k,b_k\right[$;
\item there exists a function $q_k\in C^0([a_k,b_k])\cap C^2(\left]a_k,b_k\right])$,
with $q_k(x)>0$, $q_k'(x)\geq0$ and $q_k''(x)\leq0$ on $\left]a_k,b_k\right]$,
$q_k(a_k)=0$ and $(1/q_k)\in L^1(a_k,b_k)$, and there exist two
constants $\alpha_k$ and $\beta_k$, $0<\alpha_k\leq\beta_k$, such
that $\alpha_k q_k(x)\leq p_k(x)\leq \beta_k q_k(x)$ in
$\left]a_k,b_k\right[$;
\item there exists a function $q_k\in C^0([a_k,b_k])\cap C^2(\left[a_k,b_k\right[)$,
with $q_k(x)>0$, $q_k'(x)\leq0$ and $q_k''(x)\leq0$ on $\left[a_k,b_k\right[$,
$q_k(b_k)=0$ and $(1/q_k)\in L^1(a_k,b_k)$, and there exist two
constants $\alpha_k$ and $\beta_k$, $0<\alpha_k\leq\beta_k$, such
that $\alpha_k q_k(x)\leq p_k(x)\leq \beta_k q_k(x)$ in
$\left]a_k,b_k\right[$.
\endroster
\remark{Remark} The technical condition (C)  is general enough to cover
 all the classes of nicely decomposable domains
discussed in the remarks following Definition \rft de...\endremark
\proclaim{Proposition \dft asym1..} Assume $\Omega$ satisfies
condition (C) above. Then $$
\limsup_{\nu\to\infty}{{\lambda^\oplus_\nu}\over{\nu^2}}>0. $$
\endproclaim
\demo{Proof} Condition (C), together with a simple renormalization
of the interval $\left]a_k,b_k\right[$ and a straightforward application of
Proposition \rft t2.1.., imply that for every $k=1$, \dots, $r$ we
can find a constant $\gamma_k>0$ such that
$\lambda^k_\nu\geq\gamma_k (\nu-1)^2$ for $\nu=1,2,\dots$ Let
$\gamma:=\min\{\gamma_1,\dots,\gamma_k\}$. Then, for $k=1$, \dots,
$r$ and $\nu\in\N$, we have $m(k,\gamma(\nu-1)^2)\leq \nu-1$. It follows
that $m(\oplus,\gamma(\nu-1)^2)=\sum_{k=1}^rm(k,\gamma(\nu-1)^2)\leq
r(\nu-1)$ for all $\nu\in\N$. But this means that
$\lambda^\oplus_{r(\nu-1)+1}\geq \gamma(\nu-1)^2$ for all $\nu\in\N$.
Accordingly,
we have: $$
\limsup_{\nu\to\infty}{{\lambda^\oplus_{r(\nu-1)+1}}\over{(r(\nu-1)+1)^2}}
\geq\limsup_{\nu\to\infty}{{\gamma(\nu-1)^2}\over{(r(\nu-1)+1)^2}}
={\gamma\over{r^2}}>0 $$ and the conclusion follows.\qed
\enddemo
\proclaim{Proposition \dft myp0..}
Let $(\mu_\nu)_{\nu\in\N}$ be a nondecreasing sequence of nonnegative real
numbers such that
$$\limsup_{\nu\to \infty}\frac {\mu_\nu}{\nu^2}>0 $$ and
$$\mu_\nu\le \beta^2\nu^2 \quad \text{for some $\beta\in \left]0,\infty\right[$
 and all
$\nu\in\N$.}$$
Then
$$\limsup_{\nu\to \infty}\frac{
\mu_{\nu+1}-\mu_{\nu}}{\mu_{\nu}{}^{1/2}}>0.
$$

\endproclaim
\demo{Proof} Our hypotheses imply in particular, that $\mu_\nu>0$ for all $\nu$ large
enough. Thus,              if the proposition is not true, then
$$\lim_{\nu\to \infty}\frac{
\mu_{\nu+1}-\mu_{\nu}}{\mu_{\nu}{}^{1/2}}=0.
$$
 Let $\eps>0$ be arbitrary. Then there is a $\nu_0\in\N$ such that
$$0\le\frac{\mu_{\nu+1}-\mu_{\nu}}{\mu_{\nu}{}^{1/2}}<\eps/\beta,\quad\nu\ge\nu_0.$$
Therefore
$$\mu_{\nu+1}\le\mu_{\nu}+(\eps/\beta)\mu_{\nu}{}^{1/2}\le\mu_{\nu}
+\eps \nu, \quad \nu\ge\nu_0. $$
Hence,
$$\mu_{\nu+1}\le\mu_{\nu_0}+\eps(\sum_{j=\nu_0}^\nu j)
= \mu_{\nu_0}+(\eps/2)\left(\nu(\nu+1)-\nu_0(\nu_0+1)\right),
\quad \nu\ge\nu_0.\tag \dff myfor-4.. $$
Letting $\nu\to \infty$ in \rff myfor-4.. implies that
$$\limsup_{\nu\to \infty}\frac {\mu_{\nu}}{\nu^2}\le\frac
\eps 2.$$
Since $\eps>0$ is arbitrary, it follows that
$$\limsup_{\nu\to \infty}\frac {\mu_{\nu}}{\nu^2}=0,$$ a
contradiction, which proves the proposition. \qed

       \enddemo

We can now state the main result of this section:
\proclaim{Theorem \dft gap..} Let $\Omega\subset\R^2$ be a nicely
decomposed domain and assume condition (C) is satisfied. Let
$(\lambda_\nu)_{\nu\in\N}$ be the repeated sequence of eigenvalues of
the pair $(a,b)$ in $H^1_s(\Omega)$. Then the gap condition \rff
myfk0.. is satisfied.
\endproclaim
\demo{Proof} Using Propositions~ \rft mypr-1.. and \rft mypr-2.. together with
Proposition~\rft asym1.. we see that
$$\limsup_{\nu\to \infty}\frac {\lambda_\nu}{\nu^2}>0 $$ and
$$\lambda_\nu\le \beta^2\nu^2 \quad \text{for some $\beta\in \left]0,\infty\right[$
 and all
$\nu\in\N$.}$$ Therefore Proposition~\rft myp0.. completes the proof.\qed\enddemo

\newsection\head \the\secnumber.
Inertial manifolds\endhead

Let $f\co\R\to \R$ be a function satisfying the hypotheses of
Corollary~\rfa t1.2.. together with the following dissipativeness
condition:
$$\limsup_{|s|\to \infty}f(s)/s\le -\delta_0 \quad \text{for some
$\delta_0>0$.}$$
Then, by results of \cite{\rfa pr..}, for every $\eps\ge0$ the
semiflow $\pi_\eps:=\pi_{\eps,\hat f}$ possesses a global attractor $\Cal
A_\eps$. Moreover, the family $(\Cal
A_\eps)_{\eps\ge 0}$ of attractors is upper semicontinuous at $\eps=0$.

If the eigenvalues of the limit operator $A_0$      satisfy the gap condition
\rff myfk1.. below, then, as we shall prove in Theorem~\rft myt1.. below,
 there exists an $\eps_0>0$ and a family
 $\Cal M_\eps$, $0\le \eps\le \eps_0$ of
 $C^1$-manifolds of some finite dimension $\nu$  such
 that, whenever $0\le \eps\le \eps_0$, then   $\Cal A_\eps\subset \Cal M_\eps$
 and the manifold $\Cal M_\eps$ is locally invariant relative to the
 semiflow $\pi_\eps$ on a $H^1$-neighborhood of the attractor $\Cal
 A_\eps$. Moreover, the flows on the center manifolds $\Cal M_\eps$
 converge in the (regular) $C^1$-sense to the flow on $\Cal M_0$.

 In particular, this result is valid on nicely decomposable domains satisfying condition
 $(C)$ and so this      extends and improves the corresponding results of Hale and
 Raugel \cite{\rfa HaRau1..}
 obtained for domains representable as ordinate sets with respect
 to some positive function. In fact, as it was explained in the Introduction,
  the inertial manifolds
 constructed in \cite{\rfa HaRau1..} are relative to some modified
 semiflows, which are equal to the original semiflows $\pi_\eps$ on the
 attractor $\Cal A_\eps$ but are different from $\pi_\eps$ on every
 neighborhood of it.

 The proof of Theorem~\rft myt1.. is based on the method
of functions of exponential growth, used before by many researchers (cf.
\cite{\rfa chow..}, \cite{\rfa ryb..} and the references contained in these
papers). First we choose an open set $U$ in $H^1(\om)$ which includes all the
attractors $\Cal A_\eps$, $\eps\in[0,\eps_0]$, $\eps_0>0$ small.
Then we modify the Nemitski operator $\hat f$ by finding a globally
Lipschitzian map $g\co H^1(\om)\to L^2(\om)$ with $\hat f(u)=g(u)$ for
$u\in U$.    For fixed $\eps\in[0,\eps_0]$, we seek an invariant manifold
$\Cal M_\eps$ for the modified semiflow $\pi_{\eps, g}$ in the form
$\Cal M_\eps=\Lambda_\eps(\R^\nu)$, where $\Lambda_\eps\co \R^\nu\to H^1(\om)$
 is
 a map     obtained from the contraction mapping
principle applied to a properly defined nonlinear operator $\Gamma_\eps$ defined
      on a
 certain space of maps $y\co \left]-\infty,0\right]\to H^1(\om)$
 of exponential growth. If
 the operator $\Gamma_\eps$ is a contraction then the map $\Lambda_\eps$ is
  well-defined
and $\Cal A_\eps\subset \Cal M_\eps$. It follows that $\Cal M_\eps$ is
 invariant with respect to solutions of the original semiflow $\pi_{\eps,\hat f
 }$ as long as these solutions stay in the open set $U$. One can even find an
 open set $V\subset \R^\nu$ such that for $\eps\in[0,\eps_0]$ the set $\Lambda
 _\eps(V)$ is positively invariant with respect to $\pi_{\eps,\hat f
 }$ and $\Cal A_\eps\subset \Lambda_\eps(V)\subset U$.

 The only problem is that, under the usual norm $|\cdot|_\eps$ on
 $H^1(\om)$, the operator $\Gamma_\eps$ is {\bf not} a contraction. At
 this point we use an      ingenious idea due to
Brunovsk\'y and Tere\v s\v c\'ak (see Theorem~4.1 in \cite{\rfa
brute..} and its proof) and, given positive numbers $l$ and $L$
 introduce an equivalent norm
$$\|u\|_\eps=L|u|_{L^2}+l|u|_{\eps}$$
on $H^1(\om)$. Similarly as in \cite{\rfa
brute..} we seek to choose the constants $l$ and $L$ in such a way
that the operator $\Gamma_\eps$ is a uniform contraction with
respect to the norm $\|\cdot\|_\eps$. That this is possible
follows, on the one hand, from the linear estimates contained in
Lemma~\rft mylk1.. below and, on the other hand, from the
following $C^1$-cut-off-result for Nemitski operators:

\proclaim{Theorem~\dft t1.2k..} Let $\Omega$ be an open, bounded
subset of $\R^2$, with Lipschitz boundary. Let $f\in C^1(\R\to\R)$
satisfy the growth estimate $$ |f'(s)|\leq
C(1+|s|^\beta)\quad\text{for $s\in\R$} $$ where $C$ and
$\beta\in\left[0,\infty\right[$ are arbitrary real constants.

Let  $l$ be an arbitrary positive real number and $B$ be an
arbitrary bounded subset of $H^1(\Omega)$. Then there exists an
open set $U=U(l,B)\subset H^1(\om)$ including $B$,
 a positive real number $L=L(l,B)$
and a map $g=g(l,B)\in C^1(H^1(\om)\to L^2(\om))$ with $\hat f(u)=g(u)$
for $u\in U$ and
 such that $g$ maps $H^1_s(\om)$ into $ L^2_s(\om)$ and satisfies the
 estimates
 $$\sup_{u\in H^1(\om)}|g(u)|_{L^2}<\infty\tag \dff my2..$$
$$|g(u)-g(v)|_{L^2}\le
L|u-v|_{L^2}+l|u-v|_{H^1}\quad\text{for $u$, $v\in
H^1(\Omega)$}\tag \dff my3..$$ and
 $$|Dg(u)v|_{L^2}\le
L|v|_{L^2}+l|v|_{H^1}\quad\text{for $u$, $v\in
H^1(\Omega)$.}\tag\dff my4..$$
\endproclaim
\demo{Proof} We can assume, without loss of generality, that $\beta>0$. Take a real number
$q>2$ (so in particular $0<1-2/q<1$ and $q\ge 2(\beta+1)$). Choose $\theta$ with
$1-2/q<\theta<1$. By the Gagliardo-Nirenberg inequality there is
a constant $C_1$ such that
$$ |u|_{L^q}\le C_1|u|_{H^1}{}^\theta|u|_{L^2}{}^{1-\theta}
\quad \text{for $u\in H^1(\om)$. }$$ Hence, for every $\rho>0$
$$ |u|_{L^q}\le C_1\theta\rho^{1-\theta}|u|_{H^1}+C_1(1-\theta)
\rho^{-\theta}|u|_{L^2} \quad \text{for $u\in H^1(\om)$.
}\tag\dff my1..$$ In particular, this implies that there is a real
positive constant $M$ with
$$ |u|_{L^q}<M \quad \text{for $u\in B$.
}$$ Define $U$ to be the set of all $u\in H^1(\om)$ with $
|u|_{L^q}<M$. It follows that $U$ is open in  $H^1(\om)$ and
$B\subset U$. We now apply Theorem \rft t1.1.. and obtain that
$\hat f$  is a $C^1$-map from $L^q(\om)$ to $L^2(\om)$ such that
both $|\hat f(u)|_{L^2}$ and $|D\hat f(u)|_{L(L^q(\om)\to
L^2(\om))} $ are bounded on bounded subsets of $L^q(\om)$. It is
well-known                  that there is a bounded  function $h\in
C^1(L^q(\om)\to \R)$ having globally bounded  Fr\'echet derivative
and such that $h(u)=1$ if $|u|_{L^q}\le M$ and $h(u)=0$ if
$|u|_{L^q}>2M $. This can also easily be proved using Theorem
\rft t1.1... In fact, since the function $s\mapsto |s|^q$ is in
$C^1(\R\to\R)$ with derivative $s\mapsto q|s|^{q-1}$, Theorem
~\rft t1.1.. implies that the map $\xi\co u\mapsto|u|^q$ is in
$C^1(L^q(\om)\to L^1(\om))$. Furthermore, on bounded sets in $L^q(\om)$ the
map $\xi$ is
bounded and has bounded Fr\'echet derivative. Let $\phi\in
C^1(\R\to\R)$ be such that $\phi(x)=1$ if $|x|\le M^{q}$ and
$\phi(x)=0$ if $|x|>(2M)^{q}$. It follows that the function $h(u):=\phi
(|u|_{L^q}{}^q)$ for $u\in L^q(\om)$ has the desired
properties. Define $g(u)=h(u)\hat f(u)$ for $u\in L^q(\om)$. It
follows using Leibniz rule that $g$ is a $C^1$-map from $L^q(\om)$
to $L^2(\om)$ such that both $| g(u)|_{L^2}$ and
$|Dg(u)|_{L(L^q(\om)\to L^2(\om))} $ are globally bounded. It
follows that $g$ is globally Lipschitzian from $L^q(\om)$ to
$L^2(\om)$ with a Lipschitz constant $K$. Using \rff my1.. it
follows that for every $\rho>0$ and all $u$, $v\in H^1(\Omega)$
$$|g(u)-g(v)|_{L^2}\le
C_1K(1-\theta) \rho^{-\theta}|u-v|_{L^2}+
C_1K\theta\rho^{1-\theta}|u-v|_{H^1}.$$ Choose $\rho>0$ so
small that $C_1K\theta\rho^{1-\theta}\le l$ and set
$L=C_1K(1-\theta) \rho^{-\theta}$.                  Since $H^1(\om)$ is
continuously contained in $L^q(\om)$ it follows that $g$ is of
class $C^1$ as a map from $H^1(\om)$ to $L^2(\om)$ and that \rff
my2.. and \rff my3.. are satisfied. Moreover, $g$
 maps $H^1_s(\om)$ into $L^2_s(\om)$ since $\hat f$ does so.
  The estimate \rff my3..
easily implies \rff my4... Finally, the definition of $U$
immediately implies that  $\hat f(u)=g(u)$ for $u\in U$. The
theorem is proved.
 \qed\enddemo

In order to state our main result,
we need some notation:

For every $\eps\in[0,1]$ and every
$\nu\in \N$ let
$X_{\eps,\nu,1}$ be the span of the first $\nu$ eigenvectors
$w_{\eps,j}$, $j=1$, \dots, $\nu$,
 of $A_\eps$. Let
$X_{\eps,\nu,2}$ be the orthogonal complement of  $X_{\eps,\nu,1}$ in $L^2(\om)$ if
$\eps>0$ and in $L^2_s(\om)$ if $\eps=0$. Let $A_{\eps,\nu,i}$ be the
restriction of $A_\eps$ to $X_{\eps,\nu,i}$ for $i=1$, $2$.
 Let
$E_{\eps,\nu}\xi:=\sum_{j=1}^\nu\xi_jw_{\eps,j}$, $\xi\in \R^\nu$ and
$P_{\eps,\nu,i}$ be the orthogonal projection of $L^2(\om)$ onto
$X_{\eps,\nu,i}$, $i=1$, $2$ if $\eps>0$ and $P_{0,\nu,i}$ be the
orthogonal projection of $L^2_s(\om)$ onto $X_{0,\nu,i}$, $i=1$,
$2$.

\proclaim{Theorem \dft myt1..}
Suppose that $f\in C^1(\R\to\R)$ is dissipative in the sense that
$$\limsup_{|s|\to \infty}f(s)/s\le -\delta_0 \quad \text{for some
$\delta_0>0$.}$$ Furthermore, let $f$
satisfy the growth estimate $$ |f'(s)|\leq
C(1+|s|^\beta)\quad\text{for $s\in\R$,} $$ where $C$ and
$\beta\in\left[0,\infty\right[$ are arbitrary real constants.
Suppose the eigenvalues of $A_0$ satisfy the following gap
condition:
$$\limsup_{\nu\to \infty}
\frac{\lambda_{0,\nu+1}-\lambda_{0,\nu}}{\lambda_{0,\nu}{}^{1/2}}>0.\tag\dff
myfk1..$$
Then there are an $\eps_0>0$ and an open bounded set
$U\subset H^1(\om)$ such that for every $\eps\in \left[0,
\eps_0\right[$ the attractor $\Cal A_\eps$ of the
 semiflow $\pi_{\eps,\hat f}$
 lies in $U$.

 Furthermore, there exists a globally Lipschitzian map $g\in
C^1(H^1(\om)\to L^2(\om))$ with $g(u)=
\hat f(u) $ for $u\in U$.

Besides, there is a positive integer $\nu$ and
 for every $\eps\in \left[0,
\eps_0\right[$     there
is a map $\Lambda_\eps\in C^1(\R^\nu\to H^1(\om))$ if $\eps>0$ and
$\Lambda_0\in C^1(\R^\nu\to H^1_s(\om))$ such that
$$P_{\eps,\nu,1}\circ \Lambda_\eps=E_{\eps,\nu}\tag \dff myfk20..$$ and
 $\Lambda_\eps(\R^\nu)$ is an invariant manifold with
 respect to the semiflow $\pi_{\eps, g}$.

 Finally, there is an open set $V\subset \R^\nu$ such that
 for every $\eps\in \left[0,
\eps_0\right[$ $$\Cal A_\eps\subset\Lambda_\eps(V)\subset
 U$$ and the set $\Lambda_\eps(V)$
 is  positively invariant    with respect to the semiflow $\pi_{\eps,\hat f}$.

  The reduced equation on
 $\Lambda_\eps(\R^\nu)$ takes the form
 $$\dot \xi=v_\eps(\xi),\quad \xi\in\R^\nu,\tag \dff myfk23..$$
 where $$v_\eps\co \R^\nu\to \R^\nu, \quad \xi\mapsto-A_\eps E_{\eps,\nu}\xi+
 P_{\eps,\nu,1}g(\Lambda_\eps(\xi)).$$               Moreover, whenever $\eps_n\to 0^+$
and $\xi_n\to \xi_0$ in $\R^\nu$, then $$
|\Lambda_{\eps_n}(\xi_n)-\Lambda_0(\xi_0)|_{\eps_n}+
\sum_{j=1}^\nu|\partial_j\Lambda_{\eps_n}(\xi_n)-
\partial_j\Lambda_0(\xi_0)|_{\eps_n}\to 0\tag \dff myfk21..$$and
$$
|v_{\eps_n}(\xi_n)-v_0(\xi_0)|_{\R^\nu}+
\sum_{j=1}^\nu|\partial_jv_{\eps_n}(\xi_n)-
\partial_jv_0(\xi_0)|_{\R^\nu}\to 0.\tag \dff myfk22..$$\endproclaim

Given $\mu\in\R$, a  Banach space $(Y,|\cdot|_Y)$ and a function
$y\co \left]-\infty,0\right]\to Y$ we write
$$|y|_{\mu,|\cdot|_Y}:=
\sup_{t\in\left]-\infty,0\right]}e^{\mu t}|y(t)|_Y$$ and we
 denote by $BC^\mu(Y,|\cdot|_Y)$ the set of all continuous functions
 $y\co \left]-\infty,0\right]\to Y$ such that
$|y|_{\mu,|\cdot|_Y}<\infty$. The space $BC^\mu(Y,|\cdot|_Y)$
 is a Banach space with respect to the norm
 $y\mapsto |y|_{\mu,|\cdot|_Y}$.
 In particular, we write \roster
\item $BC^\mu(L^2(\om)):=BC^\mu(L^2(\om),|\cdot|_{L^2})$, with the
norm
 $|y|_{\mu, L^2}:=|y|_{\mu, |\cdot|_{L^2}} $,
\item $BC^\mu(L^2_s(\om)):=BC^\mu(L^2_s(\om),|\cdot|_{L^2})$, with
the
norm
 $|y|_{\mu, L^2}:=|y|_{\mu, |\cdot|_{L^2}}$,
\item
 $BC^\mu(H^1(\om),\eps):=BC^\mu(H^1(\om),|\cdot|_\eps)$ with the
norm
  $|y|_{\mu, \eps}:=|y|_{\mu, |\cdot|_{\eps}}$ for $\eps>0$,
 \item
 $BC^\mu(H^1(\om),0):=BC^\mu(H^1_s(\om),|\cdot|_{H^1})$ with the
norm
  $|y|_{\mu, 0}:=|y|_{\mu, |\cdot|_{H^1}}$.\endroster

  We need the following Lemma:
\proclaim{Lemma \dft mylk1..}\roster \item
For every $\eps\in[0,1]$ and every $\nu\in\N$
$$|e^{-A_{\eps,\nu,1}t}u|_{L^2}\le
e^{-\lambda_{\eps,\nu}t}|u|_{L^2}\quad u\in X_{\eps,\nu,1},\,t\le 0,$$
$$|e^{-A_{\eps,\nu,2}t}u|_{L^2}\le
e^{-\lambda_{\eps,\nu+1}t}|u|_{L^2}\quad u\in X_{\eps,\nu,2},\,t> 0,$$
$$|e^{-A_{\eps,\nu,1}t}u|_{\eps}\le (\lambda_{\eps,\nu}+1)^{1/2}
e^{-\lambda_{\eps,\nu}t}|u|_{L^2}\quad u\in X_{\eps,\nu,1},\,t\le 0,$$
$$|e^{-A_{\eps,\nu,2}t}u|_{\eps}\le \left((\lambda_{\eps,\nu+1}+1)^{1/2} +
C_{1/2} t^{-{1/2}}\right)
e^{-\lambda_{\eps,\nu+1}t}|u|_{L^2}\quad u\in X_{\eps,\nu,2},\,t> 0.$$
\item
Define the operator
$$\Xi_{\eps,\nu}(\xi,
y)(t)=e^{-A_{\eps,\nu,1}t}E_{\eps,\nu}\xi+K_{\eps,\nu}y(t)$$ where
$\xi\in\R^\nu$, and for $y\co \left]-\infty,0\right]\to L^2(\om)$ and $t\le 0$
$$K_{\eps,\nu}y(t)=\int_{0}^te^{-A_{\eps,\nu,1}(t-s)}P_{\eps,\nu,1}y(s)
\,\d s
 +\int_{-\infty}^te^{-A_{\eps,\nu,2}(t-s)}P_{\eps,\nu,2}y(s)\,\d s,
 \tag \dff myfk2..$$
 whenever     the right
  hand side of \rff myfk2.. makes sense.
Let $\zeta$
   with $\lambda_{\eps,\nu}<\zeta<\lambda_{\eps,\nu+1}$ be arbitrary.
 Then $\Xi_{\eps,\nu}$ maps $\R^\nu\times BC^\zeta(L^2(\om))$
into $BC^\zeta(H^1(\om),\eps)$ for $\eps>0$ and
$\Xi_{0,\nu}$ maps $\R^\nu\times BC^\zeta(L^2_s(\om))$
into $BC^\zeta(H^1(\om),0)$.
Moreover, for $\eps>0$ and $y\in
BC^\zeta(L^2(\om))$ (resp. for $\eps=0$ and $y\in
BC^\zeta(L^2_s(\om))$)
$$|K_{\eps,\nu}y|_{\zeta,L^2}\le
\left(\frac 1{\zeta-\lambda_{\eps,\nu}}+
\frac 1{\lambda_{\eps,\nu+1}-\zeta}\right)|y|_{\zeta,L^2},\tag\dff sgrunt1..$$
and
$$|K_{\eps,\nu}y|_{\zeta,\eps}\le
\left(\frac{(\lambda_{\eps,\nu}+1)^{1/2}}
{\zeta-\lambda_{\eps,\nu}}+
\frac{(\lambda_{\eps,\nu+1}+1)^{1/2}}{\lambda_{\eps,\nu+1}-\zeta}+
C_{1/2}'(\lambda_{\eps,\nu+1}-\zeta)^{-1/2}\right)|y|_{\zeta,L^2}.\tag\dff sgrunt2..$$
\item If $\lambda_{0,\nu}<\zeta<\lambda_{0,\nu+1}$,
 $\eps_n\to 0^+$, $\xi_n\to \xi_0$ in $\R^\nu$ and $y_n\to y_0$ in
   $BC^\zeta(L^2(\om))$, where $y_0\in BC^\zeta(L^2_s(\om))$, then,
 for all $n$ large enough, $\lambda_{\eps_n,\nu}<\zeta<
 \lambda_{\eps_n,\nu+1}$ and
 $$|\Xi_{\eps_n,\nu}(\xi_n,y_n)-\Xi_{0,\nu}(\xi_0,y_0)|_{\zeta,\eps_n}\to
 0\quad \text{as $n\to \infty$.}$$\endroster
 \endproclaim
 \demo{Proof of the lemma} The proof of part (1) is obvious
and well-known.
Part (2) follows from part (1) and well-known results about
semigroups generated by sectorial operators. We prove \rff sgrunt2.. for $\eps>0$,
the other cases being completely analogous. Let $y\in BC^\zeta(L^2(\Omega))$ and $t\leq 0$.
By the definition of $K_{\eps,\nu}$ and by the estimates of part
(1), we have
$$
\multline
|K_{\eps,\nu}y(t)|_\eps\leq (\lambda_{\eps,\nu}+1)^{1/2}e^{-\lambda_{\eps,\nu}t}|y|_{\zeta,L^2}
\left|\int_0^t e^{(\lambda_{\eps,\nu}-\zeta)s}  \,\d s\right|\\
+(\lambda_{\eps,\nu+1}+1)^{1/2}e^{-\lambda_{\eps,\nu+1}t}|y|_{\zeta,L^2}
\int_{-\infty}^t e^{(\lambda_{\eps,\nu+1}-\zeta)s}  \,\d s\\
+C_{1/2}e^{-\zeta t}|y|_{\zeta,L^2}\int_{-\infty}^t(t-s)^{-1/2}
e^{-(\lambda_{\eps,\nu+1}-\zeta)(t-s)}  \,\d s.
\endmultline
$$
Then \rff sgrunt2.. follows from simple integration.

Let us prove part (3). Suppose that $\lambda_{0,\nu}<\zeta<\lambda_{0,\nu+1}$,
 $\eps_n\to 0^+$, $\xi_n\to \xi_0$ in $\R^\nu$ and $y_n\to y_0$ in
   $BC^\zeta(L^2(\om))$.
 Since $\lambda_{\eps,\nu}\to
  \lambda_{0,\nu}$ and $\lambda_{\eps, \nu+1}\to \lambda_{0,\nu+1}$ as
  $\eps\to 0^+$ (see \cite{\rfa pr.., Th. 3.3}), it follows that there is an $n_0\in\N$ such that
    $\lambda_{\eps_n,\nu}<\zeta<
 \lambda_{\eps_n,\nu+1}$ for all $n\ge n_0$. In the rest of this proof let
 $n\ge n_0$. Set $z_n(t):=e^{\zeta t}y_n(t)$ and $z_0(t):=e^{\zeta
t}y_0(t)$ for $t\le 0$. It easily follows that
$$\sup_{r\in \left]-\infty,0\right]}|z_n(r)-z_0(r)|_{L^2}\to 0.$$
 Thus, by Theorem 3.3 in \cite{\rfa pr..}, for every $j\in\N$
$$\sup_{r\in \left]-\infty,0\right]}
 |\langle z_n(r),w_{\eps_n,j}\rangle w_{\eps_n,j} -
 \langle z_0(r),w_{0,j}\rangle w_{0,j}|_{L^2}$$ $$\le
 \sup_{r\in \left]-\infty,0\right]}|z_n(r)-z_0(r)|_{L^2}$$
 $$+2\sup_{r\in \left]-\infty,0\right]}|z_0(r)|_{L^2}|
 w_{\eps_n,j}-w_{0,j}|_{L^2}\to 0.$$
Therefore
$$\sup_{r\in \left]-\infty,0\right]}|P_{\eps_n,\nu,1}z_n(r)-
P_{0,\nu,1}z_0(r)|_{L^2}\to 0\tag \dff myfk18..$$ and so
$$\sup_{r\in \left]-\infty,0\right]}|P_{\eps_n,\nu,2}z_n(r)-
P_{0,\nu,2}z_0(r)|_{L^2}\to 0.\tag \dff myfk17..$$
Note that whenever $(a_n)$ is a sequence of positive
 numbers converging to some positive number $a$, then
 $$\sup_{r\in \left]-\infty,0\right]}|e^{a_nt}-e^{at}|\to0.\tag \dff
  myfk19..$$ Consequently, again by Theorem 3.3 in \cite{\rfa pr..},
  $$\sup_{t\in \left]-\infty,0\right]}
  e^{\zeta t}\left| e^{-A_{\eps_n
  ,\nu,1}t}E_{\eps_n,\nu}\xi_n-
  e^{-A_{0
  ,\nu,1}t}E_{0,\nu}\xi_0 \right|_{\eps_n}$$
$$\le \sum_{j=1}^\nu\sup_{t\in \left]-\infty,0\right]}\left|
e^{(\zeta-\lambda_{\eps_n,j}) t}\xi_{n,j}w_{\eps_n,j}-
e^{(\zeta-\lambda_{0,j}) t}\xi_{0,j}w_{0,j}\right|_{\eps_n}$$
$$\le \sup_{t\in \left]-\infty,0\right]}\left|
e^{(\zeta-\lambda_{\eps_n,j}) t}-e^{(\zeta-\lambda_{0,j})
t}\right|
|\xi_{0,j}w_{0,j}|_{\eps_n}$$ $$+
 \sup_{t\in \left]-\infty,0\right]}
e^{(\zeta-\lambda_{0,j}) t}
|\xi_{n,j}w_{\eps_n,j}-\xi_{0,j}w_{0,j}|_{\eps_n}\to0.$$
Furthermore,
$$\sup_{t\le0}e^{\zeta t}\left|\int_{0}^te^{-A_{\eps_n,\nu,1}
(t-s)}P_{\eps_n,\nu,1}y_n(s)\,\d s-\int_{0}^t
-e^{-A_{0,\nu,1}
(t-s)}P_{0,\nu,1}y_0(s)\,\d s\right|_{\eps_n}\le$$
$$\sum_{j=1}^\nu\sup_{t\le0}\int_{t}^0
\left|
e^{(\zeta-\lambda_{\eps_n,j})(t-s)}
\langle z_n(s),w_{\eps_n,j}\rangle w_{\eps_n,j}-
e^{(\zeta-\lambda_{0,j}) (t-s)}
\langle z_0(s),w_{0,j}\rangle w_{0,j}\right|_{\eps_n}\,\d s$$
$$\le\sum_{j=1}^\nu\sup_{t\le0}\int_{t}^0
\left|
e^{(\zeta-\lambda_{\eps_n,j})(t-s)}-
e^{(\zeta-\lambda_{0,j}) (t-s)}\right|
|\langle z_n(s),w_{\eps_n,j}\rangle w_{\eps_n,j}|_{\eps_n}\,\d s$$
$$+\sum_{j=1}^\nu\sup_{t\le0}\int_{t}^0
e^{(\zeta-\lambda_{0,j}) (t-s)}
|\langle z_n(s),w_{\eps_n,j}\rangle w_{\eps_n,j}
-\langle z_0(s),w_{0,j}\rangle w_{0,j}|_{\eps_n}\,\d s$$
$$\le \sum_{j=1}^\nu\sup_{t\le0}\int_{t}^0
\left|
e^{(\zeta-\lambda_{\eps_n,j})(t-s)}-
e^{(\zeta-\lambda_{0,j}) (t-s)}\right|\,\d s\cdot\sup_{r\le0}
|\langle z_n(r),w_{\eps_n,j}\rangle w_{\eps_n,j}|_{\eps_n}$$
$$+\sum_{j=1}^\nu\sup_{t\le0}\int_{t}^0
e^{(\zeta-\lambda_{0,j}) (t-s)}\,\d s\cdot\sup_{r\le0}
|\langle z_n(r),w_{\eps_n,j}\rangle w_{\eps_n,j}
-\langle z_0(r),w_{0,j}\rangle w_{0,j}|_{\eps_n}=:T_n.$$
Since
$$\int_{t}^0
\left|
e^{(\zeta-\lambda_{\eps_n,j})(t-s)}-
e^{(\zeta-\lambda_{0,j}) (t-s)}\right|\,\d s=
\left|\int_{t}^0
e^{(\zeta-\lambda_{\eps_n,j})(t-s)}-
e^{(\zeta-\lambda_{0,j}) (t-s)}\,\d s\right|$$
$$=\left|\frac 1{\zeta-\lambda_{\eps_n,j}}(1-
e^{(\zeta-\lambda_{\eps_n,j})(t-s)})-
\frac 1{\zeta-\lambda_{0,j}}(1-
e^{(\zeta-\lambda_{0,j})(t-s)})\right|$$
we obtain from
\rff myfk19.. and Theorem 3.3 in \cite{\rfa pr..} that
$$T_n\to0.$$ We also have
$$\sup_{t\le0}e^{\zeta t}\left|\int_{-\infty}^te^{-A_{\eps_n,\nu,2}
(t-s)}P_{\eps_n,\nu,2}y_n(s)\,\d s-\int_{-\infty}^t
-e^{-A_{0,\nu,2}
(t-s)}P_{0,\nu,2}y_0(s)\,\d s\right|_{\eps_n}$$ $$\le
\sup_{t\le0}\int_{-\infty}^t\left|e^{-A_{\eps_n,\nu,2}
(t-s)}e^{\zeta (t-s)}P_{\eps_n,\nu,2}z_n(s)
-e^{-A_{0,\nu,2}
(t-s)}e^{\zeta (t-s)}P_{0,\nu,2}z_0(s)\right|_{\eps_n}\,\d s$$
$$=\sup_{t\le0}\int_0^{\infty}\left|e^{-A_{\eps_n,\nu,2}
s}e^{\zeta s}P_{\eps_n,\nu,2}z_n(t-s)
-e^{-A_{0,\nu,2}
s}e^{\zeta s}P_{0,\nu,2}z_0(t-s)\right|_{\eps_n}\,\d s.$$
Part (1) of this lemma implies that
 there is an integrable function
$h\co\left[0,\infty\right[\to\left[0,\infty\right[$ such that
 $$\left|e^{-A_{\eps_n,\nu,2}
s}e^{\zeta s}P_{\eps_n,\nu,2}z_n(t-s)
-e^{-A_{0,\nu,2}
s}e^{\zeta s}P_{0,\nu,2}z_0(t-s)\right|_{\eps_n}\le h(s)$$
for all $n\in\N$, and $ t$, $s\le0$.

Therefore the dominated convergence theorem, together with \rff
myfk17.. and Theorem~4.1 in \cite{\rfa pr..} imply that
$$\sup_{t\le0}\int_0^{\infty}\left|e^{-A_{\eps_n,\nu,2}
s}e^{\zeta s}P_{\eps_n,\nu,2}z_n(t-s)
-e^{-A_{0,\nu,2}
s}e^{\zeta s}P_{0,\nu,2}z_0(t-s)\right|_{\eps_n}\,\d s\to0.$$
Putting everything together we obtain
$$|\Xi_{\eps_n,\nu}(\xi_n,y_n)-\Xi_{0,\nu}(\xi_0,y_0)|_{\zeta,\eps_n}\to
 0\quad \text{as $n\to \infty$.}$$
The lemma is proved.\qed\enddemo

\demo{Proof of Theorem \rft myt1..} We divide the proof into several steps:

\demo{Step 1: choice of $l$ and $U$}

For $\nu\in\N$ with $\lambda_{0,\nu+1}-\lambda_{0,\nu}>0$ define
$$\eta_\nu=(\lambda_{0,\nu+1}-\lambda_{0,\nu})/5$$
and
$$I_\nu=[\lambda_{0,\nu}+2\eta_\nu,\lambda_{0,\nu}+3\eta_\nu ].$$
It follows that
$$ \zeta-\lambda_{0,\nu}>(\lambda_{0,\nu+1}-\lambda_{0,\nu})/3
\text{ and $\lambda_{0,\nu+1}-\zeta>(\lambda_{0,\nu+1}-
\lambda_{0,\nu})/3$ for $\zeta \in I_\nu$.}$$
Hence
$$\sup_{\zeta\in I_\nu}
\left(\frac 1{\zeta-\lambda_{0,\nu}}+
\frac 1{\lambda_{0,\nu+1}-\zeta}\right)<
C_{ \nu,1}:=6\frac 1{\lambda_{0,\nu+1}-\lambda_{0,\nu}}\tag \dff myfk5.. $$ and
$$\sup_{\zeta\in I_\nu}
\left(\frac{(\lambda_{0,\nu}+1)^{1/2}}
{\zeta-\lambda_{0,\nu}}+
\frac{(\lambda_{0,\nu+1}+1)^{1/2}}{\lambda_{0,\nu+1}-\zeta}+
C'_{1/2}(\lambda_{0,\nu+1}-\zeta)^{-1/2}\right)<C_{ \nu,2},
\tag \dff myfk6..$$ where $$
C_{ \nu,2}:=3\frac{(\lambda_{0,\nu}+1)^{1/2}}
{\lambda_{0,\nu+1}-\lambda_{0,\nu}}+3
\frac{(\lambda_{0,\nu+1}+1)^{1/2}}{\lambda_{0,\nu+1}-\lambda_{0,\nu}}+
3^{1/2} C'_{1/2} (\lambda_{0,\nu+1}-
\lambda_{0,\nu})^{-1/2}.
$$
In view of \rff myfk1.. there is a $C_1\in\left[0,\infty\right[$ and a
 strictly increasing sequence
$(\nu_k)_{k\in\N}$ in $\N$ such that
$\lambda_{0,\nu_k+1}-\lambda_{0,\nu_k}>0$ for $k\in\N$ and
$$\lambda_{0,\nu_k}{}^{1/2}/(\lambda_{0,\nu_k+1}-\lambda_{0,\nu_k})\to
C_1 \quad \text{as $k\to \infty$.}$$
Since $\lambda_{0,\nu}\to \infty$ as $\nu\to \infty$ it follows that
$$\frac 1{\lambda_{0,\nu_k+1}-\lambda_{0,\nu_k}}\to 0
\quad \text{as $k\to \infty$.}$$
Consequently,
$$C_{ \nu_k,1}\to 0 \quad
 \text{and $C_{ \nu_k,2}\to 6C_1$ as $k\to \infty$.}\tag \dff
 myfk3..$$
Choose $l$ such that
 $$0<l(6C_1+1)<\frac 1 4.$$ By
Theorem~5.10 in \cite{\rfa pr..} there is an $\eps_0$ with $0<\eps_0<1$
 and a bounded set $B_1$ in
$H^1(\om)$ such that for every $\eps\in[0,\eps_0]$ the attractor
$\Cal A_\eps$ lies in $B_1$. Let $V_0$ be the Liapunov function of $\pi_
{0,\hat f}$ defined by
$$V_0\co H^1_s(\om)\to \R,\quad u\mapsto (1/2)|\nabla u|_{L^2}{}^2-
\int_\om \hat F(u)\,\d x\d y.$$
Here, as usual, $\hat F\co H^1(\om)\to L^1(\om)$ is the Nemitski operator
defined by the function $F(x):=\int_0^x f(s)\,\d s$, $x\in \R$.
There is an $M_0\in \left]0,\infty\right[$ so that $V_0(u)<M_0$ for all
$u\in \Cal A_0$.
As in \cite{\rfa pr..} it follows that
$$B_2:=\{\,u\in H^1_s(\om)\mid V_0(u)\le M_0\,\}$$ is bounded.
Define $B:=B_1\cup B_2$ and
 let $L=L(l,B)$ and $U=U(l,B)$ and $g=g(l,B)$                  be as in
Corollary \rft t1.2k... By \rff myfk3..  there exists a
 $k\in \N$ such that
$$LC_{ \nu_k,1}<\frac 1 4 \quad
 \text{and $lC_{ \nu_k,2}<\frac 1 4$.}\tag \dff myfk4..$$
  Fix such a $k$ and set $\nu:=\nu_k$. Since $\lambda_{\eps,\nu}\to
  \lambda_{0,\nu}$ and $\lambda_{\eps, \nu+1}\to \lambda_{0,\nu+1}$ as
  $\eps\to 0^+$ and using \rff myfk5.. and \rff myfk6.., we may
  assume,
by taking $\eps_0$ smaller if necessary,
that for every $\eps\in[0,\eps_0]$ and every $\zeta\in
 I_\nu$
 $$ \lambda_{\eps,
 \nu}<\zeta<\lambda_{\eps,\nu+1}
 ,$$
 $$\sup_{\zeta\in I_\nu}
\left(\frac 1{\zeta-\lambda_{\eps,\nu}}+
\frac 1{\lambda_{\eps,\nu+1}-\zeta}\right)<
C_{ \nu,1}\tag \dff myfk7.. $$ and
$$\sup_{\zeta\in I_\nu}
\left(\frac{(\lambda_{\eps,\nu}+1)^{1/2}}
{\zeta-\lambda_{\eps,\nu}}+
\frac{(\lambda_{\eps,\nu+1}+1)^{1/2}}{\lambda_{\eps,\nu+1}-\zeta}+
C'_{1/2}(\lambda_{\eps,\nu+1}-\zeta)^{-1/2}\right)<
C_{\nu,2}.
\tag \dff myfk8..$$
\enddemo

\demo{Step 2: choice of a new norm}

If $\eps>0$ endow $H^1(\om)$ with the equivalent norm
$$\|u\|_\eps:=L|u|_{L^2}+l|u|_\eps.$$
Write
$Z^\zeta_\eps:=BC^\zeta(H^1(\om),\|\cdot\|_\eps)$ with the
corresponding norm
$$\|y\|_{\zeta,\eps}=\sup_{t\in\left]-\infty,0\right]}e^{\zeta
t}\|y(t)\|_\eps.$$
If $\eps=0$ endow $H^1_s(\om)$ with the equivalent norm
$$\|u\|_0:=L|u|_{L^2}+l|u|_{H^1}.$$
Write
$Z^\zeta_0:=BC^\zeta(H^1_s(\om),\|\cdot\|_0)$ with the
corresponding norm
$$\|y\|_{\zeta,0}=\sup_{t\in\left]-\infty,0\right]}e^{\zeta
t}\|y(t)\|_0.$$
It follows that for $\eps\in[0,\eps_0]$, $\zeta\in I_\nu$ and $y\in
BC^\zeta(L^2(\om))$ (resp. for $\eps=0$ and $y\in
BC^\zeta(L^2_s(\om))$)
$$\|K_{\eps,\nu}y\|_{\zeta,\eps}\le
(1/2)|y|_{\zeta,L^2}.\tag \dff myfk15..$$
\enddemo

\demo{Step 3: the contraction map $\Gamma_\eps$}

Fix $\zeta$ and $\mu\in I_\nu$ with $\zeta<\mu$.
Since $g\co H^1(\om)\to L^2(\om)$ is globally bounded, it follows
that $g\circ y$ is globally bounded for every
 $y\in BC^\zeta(H^1(\om),\eps)$, $\eps\ge 0$. Moreover, since $g$
 maps $H^1_s(\om)$
into $L^2_s(\om)$ it follows that the nonlinear operator
 $y\mapsto g\circ y$ maps $Z^\zeta_\eps$ into
 $ BC^\zeta(L^2(\om))$ for $\eps>0$ and it maps
 $ Z^\zeta_0$ into
 $ BC^\zeta(L^2_s(\om))$.
 Moreover,
 $$|g(u)-g(v)|_{L^2}\le L|u-v|_{L^2}+l|u-v|_{H^1}$$ $$\le \|u-v\|_{\eps}\quad
 \text{for $u$, $v\in H^1(\om)$}$$ so
 $$|g\circ y-g\circ w|_{\zeta,L^2}\le \|y-w\|_{\zeta,\eps}\quad
 \text{for $y$, $w\in Z^\zeta_\eps$.}\tag \dff myfk16..$$
 It follows that the operator
 $$\Gamma_\eps\co \R^\nu\times Z^\zeta_\eps\to Z^\zeta_\eps,\quad
 (\xi,y)\mapsto \Xi_{\eps,\nu}(\xi,g\circ y)$$ is well-defined.
 Spelled out in full, the operator $\Gamma_\eps$ reads as
 follows:
 $$ \Gamma_\eps(\xi, y)(t)=e^{-A_{\eps,\nu,1}t}E_{\eps,\nu}\xi+
 \int_{0}^t e^{-A_{\eps,\nu,1}(t-s)}P_{\eps,\nu,1}g(y(s))\,\d s$$
 $$+\int_{-\infty}^te^{-A_{\eps,\nu,2}(t-s)}P_{\eps,\nu,2}g(y(s))\,\d s\quad
 \text{for $(\xi,y)\in                  \R^\nu\times Z^\zeta_\eps$ and $t\le 0$.
 }$$
If $0\le \eps\le\eps_0$ then, by \rff myfk15.. and
 \rff myfk16..,
$\Gamma_\eps$ is a uniform
 contraction in the second variable with contraction constant
 $1/2$. It follows that for every such $\eps$ there is a
  uniquely defined map
 $\phi_\eps\co \R^\nu\to Z^\zeta_\eps$, such that
 $$\phi_\eps(\xi)=\Gamma_\eps(\xi,\phi_\eps(\xi))
 \quad\text{for $\xi\in\R^\nu$.}$$
 Either proceeding directly, or by using the fiber contraction
 theorem as in \cite{\rfa chow..} or else by using the abstract results of
\cite{\rfa ryb.. }
 one proves that the map $\phi_\eps$ is of class $C^1$ as a map
 from $\R^\nu$ into $Z^\mu_\eps$ and for every $\xi'\in\R^\nu$ and $j=1$, \dots, $\nu$,
 $D\phi_\eps(\xi)\xi'$ lies in $Z^\zeta_\eps$ and is given by
 the recursive formula
 $$D\phi_\eps(\xi)\xi'=\Xi_{\eps,\nu}(\xi',v)\tag \dff myfk11..$$
 where
 $$ v(s):=Dg(\phi_\eps(\xi)(s))(D\phi_\eps(\xi)\xi'(s)),\quad s\le
 0.$$
\enddemo

\demo{Step 4: behaviour as $n\to\infty$}

 Let $\xi_n\to \xi_0$ in $\R^\nu$ and let $\eps_n\to 0^+$. We
 claim that
 $$\|\phi_{\eps_n}(\xi_n)-\phi_0(\xi_0)\|_{\zeta,\eps_n}\to 0\quad
 \text{for $n\to \infty$,}\tag \dff myfk9..$$ and for every $\xi'\in\R^\nu$
  $$\|D\phi_{\eps_n}(\xi_n)(\xi')-D\phi_0(\xi_0)(\xi')\|_{\mu,\eps_n}\to 0\quad
 \text{for $n\to \infty$.}\tag\dff myfk10..$$
 In fact,
 $$\phi_{\eps_n}(\xi_n)-\phi_0(\xi_0)=\Gamma_{\eps_n}
 (\xi_n,\phi_{\eps_n}(\xi_n))-
 \Gamma_{\eps_n}
 (\xi_n,\phi_{0}(\xi_0))$$ $$+\Gamma_{\eps_n}
 (\xi_n,\phi_{0}(\xi_0))-\Gamma_{0}
 (\xi_0,\phi_{0}(\xi_0)).$$
 and
 $$\|\Gamma_{\eps_n}
 (\xi_n,\phi_{\eps_n}(\xi_n))-
 \Gamma_{\eps_n}
 (\xi_n,\phi_{0}(\xi_0))\|_{\zeta,\eps_n}\le (1/2)
 \|\phi_{\eps_n}(\xi_n)-\phi_0(\xi_0)\|_{\zeta,\eps_n}$$
 so
 $$\|\phi_{\eps_n}(\xi_n)-\phi_0(\xi_0)\|_{\zeta,\eps_n}\le
 2\|\Gamma_{\eps_n}
 (\xi_n,\phi_{0}(\xi_0))-\Gamma_{0}
 (\xi_0,\phi_{0}(\xi_0))\|_{\zeta,\eps_n}$$
 $$=2\|\Xi_{\eps_n,\nu}(\xi_n,y_n)-\Xi_{0,\nu}(\xi_0,y_0)\|_{\zeta,\eps_n}$$
 $$\le
 2(L+l)|\Xi_{\eps_n,\nu}(\xi_n,y_n)-\Xi_{0,\nu}(\xi_0,y_0)|_{\zeta,\eps_n}$$
 where $y_n\equiv y_0=g\circ(\phi_{0}(\xi_0))$. Therefore \rff
 myfk9..
 follows from                  Lemma \rft mylk1...
 Now
 $$|Dg(u)v|_{L^2}\le L|v|_{L^2}+l|v|_{H^1}\le \|v\|_{\eps}\quad
 \text{for $u$, $v\in H^1(\om)$.}\tag \dff myfk14..$$
 Therefore, using \rff myfk11..  and proceeding as in the proof
 of \rff myfk9.., we obtain
$$\gather\|D\phi_{\eps_n}(\xi_n)(\xi')-D\phi_0(\xi_0)(\xi')\|_{\mu,\eps_n}\\
\le2\|\Xi_{\eps_n,\nu}(\xi',y_n)-\Xi_{0,\nu}(\xi',y_0)\|_{\mu,\eps_n}\\
\le 2(L+l)|\Xi_{\eps_n,\nu}(\xi',y_n)-\Xi_{0,\nu}(\xi',y_0)
|_{\mu,\eps_n}
\tag \dff myfk13..\endgather$$ where,
for every $s\le 0$,
$$v(s):=D\phi_0(\xi_0)\xi'(s),$$
$$y_n(s):=Dg(\phi_{\eps_n}(\xi_n)(s))v(s)$$
and
$$y_0(s):=Dg(\phi_{0}(\xi_0)(s))v(s).$$
We claim that
$$|y_n-y_0|_{\zeta,L^2}\quad\text{as $n\to \infty$.}
\tag \dff myfk12..$$
Assuming this claim, we obtain \rff myfk10.. from \rff myfk13..,
\rff myfk12.. and                  Lemma \rft mylk1...
To prove the claim, note that by \rff myfk14..
$$e^{\mu t}|y_n(t)-y_0(t)|_{L^2}=e^{(\mu-\zeta) t}
|Dg(\phi_{\eps_n}(\xi_n)(t))w(t)-Dg(\phi_{0}(\xi_0)(t))w(t)|_{L^2}$$
$$\le e^{(\mu-\zeta)t}\|v\|_{\eta,0}$$
where $w(t):=e^{\eta t}v(t)$, so $\|w(t)\|_{\eps_n}
\|w(t)\|_{0}\le \|v\|_{\eta,0}$,
  $t\le 0$.

Therefore, given $\delta>0$ there is a $a\le 0$ such that
$$e^{\mu t}|y_n(t)-y_0(t)|_{\eps_n}\le \delta
\quad \text{for $t<a$.}$$
Thus we only need to prove that
$$\sup_{t\in[a,0]}|Dg(\phi_{\eps_n}
(\xi_n)(t))w(t)-Dg(\phi_{0}(\xi_0)(t))w(t)|_{L^2}\to 0
\quad\text{as $n\to \infty$.}$$
However, this is obvious since $Dg\co H^1(\om)\to
L(H^1(\om),L^2(\om))$ is continuous and, by \rff myfk9..,
$$\sup_{t\in[a,0]}|\phi_{\eps_n}
(\xi_n)(t)-\phi_{0}(\xi_0)(t)|_{H^1}\to 0
\quad\text{as $n\to \infty$.}$$
\enddemo

\demo{Step 5: the global invariant manifold}

Define, for $0<\eps\le \eps_0$,
 $$\Lambda_\eps\co \R^\nu\to H^1(\om),\quad \xi\mapsto
 \phi_\eps(\xi)(0),$$ and
 $$\Lambda_0\co \R^\nu\to H^1_s(\om),\quad \xi\mapsto
 \phi_0(\xi)(0).$$
 By what we have proved so far, for every $\eps\in[0,\eps_0]$ the map
  $\Lambda_\eps$ is well-defined, of
 class $C^1$ and  \rff myfk20.. and \rff myfk21.. hold.
 It is well-known and easily proved
  that $\Lambda_\eps(\R^\nu)$ is an invariant manifold
 of the semiflow $\pi_{\eps,g}$ which includes all orbits of
 solutions of $\pi_{\eps,g}$ defined for $t\le 0$ and lying in
 $Z^\zeta_\eps$. Since $g$ equals $\hat f$ on $U$, it follows that
  every point in
 $\Cal A_\eps$ is contained in $\Lambda_\eps(\R^\nu)$.
  The reduced                  equation on the manifold $\Lambda_\eps(\R^\nu)$ clearly
  takes the form \rff myfk23.. and \rff myfk21.. implies \rff
  myfk22...
\enddemo

\demo{Step 6: the local invariant manifold}

  Let
  $$K:=\{\,\xi\in\R^\nu\mid V_0(\Lambda_0(\xi))\le M_0\,\}=
  \{\,\xi\in\R^\nu\mid \Lambda_0(\xi)\in B_2\,\}.$$
  Since $B_2$ is bounded and closed, it follows from \rff myfk20..
  that $K$ is bounded and closed, i.e. compact.

  Define
  $$V:=\{\,\xi\in\R^\nu\mid V_0(\Lambda_0(\xi))< M_0\,\}.$$
 Thus $V\subset K$ and $V$ is open in $\R^\nu$.
 Since $\Lambda(K)\subset U$ and $K$ is compact and $U$ is open in
 $H^1(\om)$, it follows
 from \rff myfk21.., by choosing $\eps_0>0$ smaller, if necessary, that
 $$\Lambda_\eps(K)\subset U,\quad \eps\in[0,\eps_0].$$
 We also claim that, if $\eps_0>0$ is small enough, then
 $$ \Lambda_\eps{}^{-1}(\Cal A_\eps)\subset V, \quad \eps\in [0,\eps_0].$$
  In fact if the claim does not hold, then
 there are                  sequences $\eps_n\to0^+$ and $\xi_n\notin V$ with
 $u_n:=\Lambda_{\eps_n}(\xi_n)\in\Cal A_{\eps_n}$. By Corollary~5.2 and Lemma~5.9
 of \cite{\rfa pr..} we may assume that $u_n\to u_0$ in $H^1(\om)$, where
 $u_0\in \Cal A_0$. It follows that there is a $\xi_0\in \R^\nu$ with
 $u_0=\Lambda_0(\xi_0)$. Since, by \rff myfk20.., for every $a\in \R^\nu$ and every
 $\eps\in[0,\eps_0]$, the
 $j$-th component $a_j$ of $a$ is given by $a_j=\langle
 \Lambda_\eps(a),w_{\eps,j}\rangle$, it follows that
 $\xi_n\to \xi_0$ in $\R^\nu$. By the definition of $V$ we have
 $\xi\in V$, a contradiction, proving the claim.
Set $W:=V_0\circ\Lambda_0\co \R^\nu\to \R$.
Then, for every $\xi\in \Lambda_\eps{}^{-1}(U)$, we have
$$\nabla W(\xi)\cdot
v_0(\xi)=DV_0(\Lambda_0(\xi))D\Lambda_0(\xi)(v_0(\xi))$$
$$=-|D\Lambda_0(\xi)(v_0(\xi))|_{L^2}{}^2.$$  Since there are no
equilibria $u$ of $\pi_{0,\hat f}$ with $V_0(u)=M_0$, it follows that
 $$D\Lambda_0(\xi)(v_0(\xi))\not=0,\quad \text{whenever $W(\xi)=M_0$.}$$
 By the compactness of $K$ we now obtain
that there is a $\delta>0$ such that
$$\nabla W(\xi)\cdot
v_0(\xi)<-\delta,\quad \text{whenever $W(\xi)=M_0$.}$$
Therefore \rff myfk22.. implies that, if $\eps_0>0$ is small enough, then
$$\nabla W(\xi)\cdot
v_\eps(\xi)<-\delta,\quad \text{whenever $W(\xi)=M_0$.}$$
This shows, that, for $\eps\in[0,\eps_0]$, the set $V$ is positively invariant for the
equation \rff myfk23.. so the set
 $\Lambda_\eps(V)\subset U$ is positively invariant for the semiflow
$\pi_{\eps,\hat f}$ and $\Cal A_\eps\subset \Lambda_\eps(V)$.
The theorem is proved. \qed
\enddemo
 \enddemo

\newsection\head \the\secnumber.
Final remarks\endhead

It is clear that Theorem \rft myt1.., by its abstract nature, can  easily be
 generalized to the case of an arbitrary
smooth bounded domain
$\Omega\subset\R^M\times\R^N$, with
$n:=M+N\geq 2$.  In this case the squeezed domain $\Omega_\eps$ is defined by
$$
\Omega_\eps:=\{\,(x,\eps y)\mid (x,y)\in\Omega\,\},
$$
where now $x\in\R^M$ and $y\in\R^N$.
If $n>2$ one
only needs to impose the condition
$$\beta\le (2^*/2)-1, \quad\text{ where  $2^*=2n/(n-2)$}\tag\dff gasp..$$
in \rff f1.1.., so the
Nemitski operator associated with $f$ is a $C^1$-map of $H^1$ into $L^2$ (see \cite{\rfa pr..} for
details). The strict inequality in \rff gasp.. makes
 it possible to use the Gagliardo-Nirenberg
inequality in the proof of Theorem \rft t1.2k...

\vbox to 2in{\vss\includegraphics{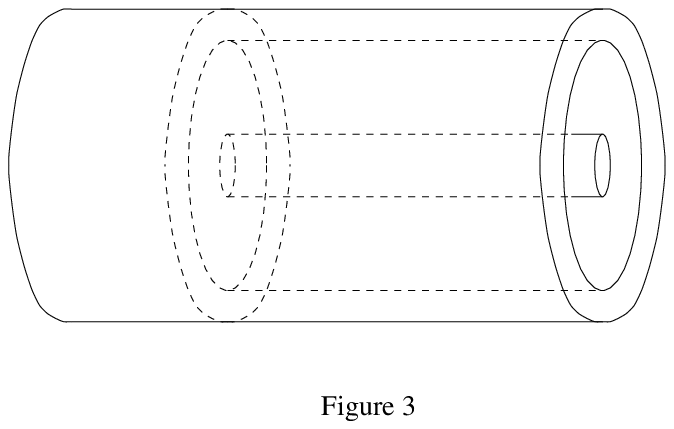}}

Of course, this generalization is meaningful only if one knows
that the eigenvalues of the limit operator $A_0$ satisfy the gap condition \rff myfk1...
Under fairly general conditions, this is actually the case if $M=1$ and
$N\geq 1$ is arbitrary. In fact, one can easily extend the notion of {\it nicely decomposed domain}
in order to cover this more general situation. If the additional condition (C) introduced in
Section 4 is satisfied, then one can prove \rff myfk1.. by just repeating the arguments
of Section 4. We do not enter into detail: Figure 3 illustrates a simple example
of a three-dimensional domain on which this kind of analysis is possible.

If $M\geq 2$ the problem of the spectral gap condition is much more delicate (see  \cite{\rfa
MaSe..} and \cite{\rfa MaSeSh..} for this and some related questions). Some (partial) generalizations
of Theorem \rft myt1.. should be possible for very particular three-dimensional domains, essentially
for domains which are ordinate sets of positive functions defined on rectangles.

One can also consider, much more generally, a squeezing transformation of an open subset
of $\R^{M+N}$ toward an arbitrary smooth $M$-dimensional submanifold $S$ of $\R^{M+N}$ (see
\cite{\rfa prr..} for details). If the eigenvalues of the Laplace-Beltrami operator on $S$
satisfy a gap condition like \rff myfk1.. (e.g. when $S$ is a sphere), some further
generalizations of Theorem \rft myt1.. have been obtained in \cite{\rfa pr4..}.

\newsection\head\the\secnumber.
Appendix
\endhead
In this appendix we give the \demo{Proof of Lemma \rft t2.2..}
First of all, we observe that, since $q''(x)\leq0$ on $\left]0,1\right]$, the
function $q$ is concave on $\left]0,1\right[$. So, for $x\in\left]0,1\right[$ and for
all $0<\epsilon<x$, we have $q(x)+(\epsilon-x)q'(x)\geq
q(\epsilon)$. Now, since $q$ is continuous on $[0,1]$ and
$q(0)=0$, we get $$ 0\leq (q'(x)/q(x))\leq (1/x)\quad \text{for
$x\in\left]0,1\right]$.}\tag \dff f2.0.. $$

Next, we analyze some properties of eigenvectors of $(a^q,b^q)$.
Let $\lambda$ be an eigenvalue of $(a^q,b^q)$, and let $u\in V(q)$
be a corresponding eigenvector. This means that $$
\int_0^1qu'v'\,\d x=\lambda\int_0^1quv\,\d x \quad\text{for all
$v\in V(q)$.} $$ Choosing $v\in C^\infty_0(\left]0,1\right[)$ arbitrarily, we
see that $(qu')'\in L^1_{\roman{loc}}(0,1)$ and $(qu')'=-\lambda
qu$. In particular, we see that $(qu')'\in C^0([0,1])$, and hence
$qu'\in C^1([0,1])$. Moreover, taking $v$ any function in
$C^\infty([0,1])$, with $v\equiv 0$ on $[0,\cl\epsilon]$, $v\equiv
1$ on $[1-\cl\epsilon,1]$, where $0<\cl\epsilon<1-\cl\epsilon$, we
obtain for any $\epsilon$, $0<\epsilon<\cl\epsilon$: $$ \align
&\lambda\int_0^1quv\,\d x =\int_0^1qu'v'\,\d
x=\int_\epsilon^{1-\epsilon}qu'v'\,\d x\\
&=-\int_\epsilon^{1-\epsilon}(qu')'v\,\d x
+q(1-\epsilon)u'(1-\epsilon)v(1-\epsilon)-q(\epsilon)u'(\epsilon)v(\epsilon)\\
&=\lambda\int_\epsilon^{1-\epsilon}quv\,\d
x+q(1-\epsilon)u'(1-\epsilon).
\endalign
$$ By letting $\epsilon\to 0$, we obtain that $u'(1)=0$. Taking
now $v$ any function in $C^\infty([0,1])$, with $v\equiv 1$ in
$[0,\cl\epsilon]$, $v\equiv 0$ on $[1-\cl\epsilon,1]$, where
$0<\cl\epsilon<1-\cl\epsilon$, we obtain for any $\epsilon$,
$0<\epsilon<\cl\epsilon$: $$ \align &\lambda\int_0^1quv\,\d x
=\int_0^1qu'v'\,\d x=\int_\epsilon^{1-\epsilon}qu'v'\,\d x\\
&=-\int_\epsilon^{1-\epsilon}(qu')'v\,\d x
+q(1-\epsilon)u'(1-\epsilon)v(1-\epsilon)-q(\epsilon)u'(\epsilon)v(\epsilon)\\
&=\lambda\int_\epsilon^{1-\epsilon}quv\,\d
x-q(\epsilon)u'(\epsilon).
\endalign
$$ By letting $\epsilon\to 0$, we obtain that $(qu')(0)=0$. But
now, for $x\in\left]0,1\right]$, we have: $$ (qu')(x)=\int_0^x(qu')'(s)\,\d
s=-\lambda\int_0^xq(s)u(s)\,\d s; $$ since $q$ is nonnegative and
nondecreasing on $[0,1]$, it follows that $$
|q(x)u'(x)|\leq\lambda\int_0^xq(s)|u(s)|\,\d s\leq\lambda
q(x)|u|_{\infty}x, $$ and, since $q(x)>0$ for $x\in\left]0,1\right]$, $$
|u'(x)|\leq\lambda|u|_\infty x.\tag\dff f2.1.. $$ In particular,
$u\in C^1([0,1])$ and $u'(0)=0$.

Assume now that $\lambda\not=0$. Then $u'\not\equiv0$. Let us
define $$ w(x):=q(x)^{1/2}u'(x),\quad x\in[0,1]. $$ Then $w\in
C^0([0,1])$ and $w(0)=w(1)=0$. Moreover, in $\left]0,1\right]$ we have: $$
\align
w'&=(q^{-1/2}qu')'=-{{1}\over{2}}q^{-3/2}q'qu'+q^{-1/2}(qu')'\\
&=-{{1}\over{2}}{{q'}\over{q}}q^{1/2}u'-q^{-1/2}\lambda qu\\
&=-{{1}\over{2}}{{q'}\over{q}}w-\lambda q^{1/2}u.\\
\endalign
$$ From this equality, it follows that $w'\in C^0(\left]0,1\right])$, so
$w\in C^1(\left]0,1\right])$; but then the same equality implies that $w'\in
C^1(\left]0,1\right])$, so $w\in C^2(\left]0,1\right])$. Moreover, \rff f2.0.. and \rff
f2.1.. imply that $w'\in L^2(0,1)$, so $w\in H^1_0$. We compute
$w''$: $$ \align
w''&=-{{1}\over{2}}\left({{q'}\over{q}}\right)'w-{{1}\over{2}}{{q'}\over{q}}w'-
\lambda (q^{1/2}u)'\\
&=-{{1}\over{2}}\left({{q'}\over{q}}\right)'w-{{1}\over{2}}{{q'}\over{q}}
\left(-{{1}\over{2}}{{q'}\over{q}}w-\lambda q^{1/2}u\right)
-{{1}\over{2}}\lambda q^{-1/2}q'u-\lambda w\\
&=-{{1}\over{2}}\left({{q'}\over{q}}\right)'w+{{1}\over{4}}{{q'\null^2}\over{q^2}}w
-\lambda w\\ &=-{{1}\over{2}}{{q''q-q'\null^2}\over{q^2}}w
+{{1}\over{4}}{{q'\null^2}\over{q^2}}w-\lambda w\\
&=\left({{3}\over{4}}{{q'\null^2}\over{q^2}}-
{{1}\over{2}}{{q''}\over{q}}\right)w-\lambda w\\
\endalign
$$ Let us define $$
Q(x):={{3}\over{4}}{{q'(x)\null^2}\over{q(x)^2}}-
{{1}\over{2}}{{q''(x)}\over{q(x)}},\quad x\in \left]0,1\right]. $$ Observe
that $0\leq Q(x)\leq -(3/4)(q'/q)'(x)$ on $\left]0,1\right]$. We claim that
$Q^{1/2}w\in L^2(0,1)$. In fact let $\epsilon$, $0<\epsilon<1$, be
arbitrary: by \rff f2.0.. and \rff f2.1.., we have: $$ \align
\int_\epsilon^1Qw^2\,\d x
&\leq{3\over4}\int_\epsilon^1\left(-{{q'}\over{q}}\right)'w^2\,\d
x\\ &\leq{3\over4}\int_\epsilon^1\left(-{{q'}\over{q}}\right)'
(\lambda|u|_\infty|q|_\infty^{1/2} x)^2\,\d x\\
&\leq{3\over4}\lambda|u|_\infty^2|q|_\infty
\int_\epsilon^1\left(-{{q'}\over{q}}\right)'x^2\,\d x\\
&={3\over4}\lambda|u|_\infty^2|q|_\infty
\left(\left.-{{q'}\over{q}}x^2\right|^1_\epsilon+2\int_\epsilon^1
{{q'}\over{q}}x\,\d x\right)\\ &\leq
{3\over4}\lambda|u|_\infty^2|q|_\infty
\left({1\over\epsilon}\epsilon^2+2\int_\epsilon^1 {1\over x}x\,\d
x\right).\\
\endalign
$$ Since $\epsilon$ is arbitrary, the claim is proved.

Let us define the Hilbert space $$ U(Q):=\left\{\,\omega\in
H^1_0(0,1)\mid Q^{1/2}\omega\in L^2(0,1)\,\right\}, $$ with the
scalar product $$ \langle\omega,\varpi\rangle_{U(Q)}:=
\int_0^1\omega'\varpi'\,\d x+\int_0^1Q\omega\varpi\,\d x. $$ Let
us notice that $U(Q)\hookrightarrow H^1_0(0,1)$ with continuous
imbedding, so the imbedding $U(Q)\hookrightarrow L^2(0,1)$ is
compact. Moreover, $U(Q)$ is clearly a dense subspace of
$L^2(0,1)$, since it contains $C^\infty_0(\left]0,1\right[)$. The hypothesis
of Proposition~2.2 in \cite{\rfa pr..} is then satisfied by the pair $({\roman
a}_Q, {\roman b})$, where ${\roman
a}_Q:=\langle\cdot,\cdot\rangle_{U(Q)}$ and ${\roman
b}:=\langle\cdot,\cdot\rangle_{L^2(0,1)}$ . Then it follows that
there exists an eigenvalue-eigenvector sequence $(\mu_\nu^Q,\omega_\nu^Q)_{\nu\in\N}$ of
the pair $({\roman a}_Q, {\roman b})$ such
that $$ 0\leq\mu_1^Q\leq\mu_2^Q\leq\mu_3^Q\leq\dots $$ and
$(\omega^Q_\nu)_{\nu\in\N}$ is a complete orthonormal system in
$L^2(0,1)$. We have proved above that $w\in U(Q)$. We claim that
$$ \int_0^1w'\omega'\,\d x+\int_0^1Qw\omega\,\d x=
\lambda\int_0^1w\omega\,\d x\quad\text{for all $\omega\in U(Q)$,}
$$ i.e. $(\lambda, w)$ is an eigenvalue-eigenvector pair of
$({\roman a}_Q, {\roman b})$. In fact, for every $\epsilon$,
$0<\epsilon<1$, we have $$ \align \int_\epsilon^1w'\omega'\,\d
x&=-\int_\epsilon^1w''\omega\,\d x+
\left.w'\omega\right|^1_\epsilon\\ &=-\int_\epsilon^1Qw\omega\,\d
x+\lambda\int_\epsilon^1w\omega\,\d x
-w'(\epsilon)\omega(\epsilon).\\
\endalign
$$ Then $$ \int_0^1w'\omega'\,\d x+\int_0^1Qw\omega\,\d x=
\lambda\int_0^1w\omega\,\d x-\lim_{\epsilon\to
0}w'(\epsilon)\omega(\epsilon). $$ Recall that $$ w'(\epsilon)=
\left(-{1\over2}{{q'(\epsilon)}\over{q(\epsilon)}}
q(\epsilon)^{1/2}u'(\epsilon)- \lambda
q(\epsilon)^{1/2}u(\epsilon)\right); $$ so, by \rff f2.0.. and
\rff f2.1.., it follows that $w'(\epsilon)\omega(\epsilon)\to 0$
as $\epsilon\to 0$ and the claim is proved.

Thus we have proved that, whenever $(\lambda,u)$ is an
eigenvalue-eigenvector pair for $(a^q,b^q)$ and $\lambda\not=0$,
then $(\lambda,q^{1/2}u')$ is an eigenvalue-eigenvector pair for
$({\roman a}_Q, {\roman b})$. Assume now that $\lambda\not=0$ is
an eigenvalue of $(a^q,b^q)$, and that $u,\tilde u\in V(q)$ are
two linearly independent corresponding eigenvectors. Let
$w:=q^{1/2}u'$ and $\tilde w:=q^{1/2}\tilde u'$. Then also $w$ and
$\tilde w$ are linearly independent. Otherwise, we could find
$\xi\not=0$ such that $w=\xi\tilde w$, that is $q^{1/2}u'=\xi
q^{1/2}\tilde u'$. Since $q(x)>0$ on $\left]0,1\right]$, it follows that
$u'=\xi\tilde u'$. So there is a constant $\zeta$ such that $u=\xi
\tilde u+\zeta$, that is $u$, $\tilde u$ and $1$ are linearly
dependent, a contradiction. This means that the multiplicity of
$\lambda$ as an eigenvalue of $(a^q,b^q)$ is less then or equal to the
multiplicity of $\lambda$ as an eigenvalue of $({\roman a}_Q,
{\roman b})$. So we can conclude that $$
\lambda^q_1=0\quad\text{and}\quad
\lambda^q_\nu\geq\mu_{\nu-1}^Q\quad\text{for $\nu\geq 2$.} $$ In order
to complete the proof, we need to show that $\mu_\nu^Q\geq\pi^2\nu^2$,
for $\nu=1,2,\dots$ This is done as follows: by Proposition~2.2 in \cite{\rfa pr..}, $$
\mu_\nu^Q=\inf_{E\in\Cal U_\nu(Q)}\sup_{\omega\in
E\setminus\{0\}} {{\int_0^1\omega'\null^2\,\d
x+\int_0^1Q\omega^2\,\d x} \over{\int_0^1\omega^2\,\d x}}, $$
where $\Cal U_\nu(Q)$ is the set of all $\nu$-dimensional subspaces of
$U(Q)$. Since $Q\geq0$ and $U(Q)\subset H^1_0(0,1)$, it follows
that $$ \mu_\nu^Q\geq\mu_\nu^0:=\inf_{E\in\Cal U_\nu}\sup_{\omega\in
E\setminus\{0\}} {{\int_0^1\omega'\null^2\,\d x}
\over{\int_0^1\omega^2\,\d x}}, $$ where $\Cal U_\nu$ is the set of
all $\nu$-dimensional subspaces of $H^1_0(0,1)$. But
$(\mu_\nu^0)_{\nu\in\N}$ are exactly the eigenvalues of the Dirichlet
problem $$ \cases \omega''(x)=\mu\omega(x),&x\in\left]0,1\right[\\
u(x)=0,&x\in\{0,1\}
\endcases
$$ This implies that $\mu_\nu^0=\pi^2\nu^2$, $\nu=1,2,\dots$, and the
proof is complete. \qed\enddemo

\enddocument\bye